%% file: UniformGeomSpace-en.tex
\documentclass[12pt,a4paper]{book}
\usepackage{amsmath,amsfonts,amssymb}
\usepackage[utf8]{inputenc}
\usepackage{graphicx}
\usepackage{wrapfig}
\usepackage[pdftex,colorlinks]{hyperref}
\begin{document}
\begin{titlepage}
\vfill
\begin{center}
\textbf{\Huge Uniform Theory\\ of Geometric Spaces}

\vspace{1.5cm}
\includegraphics[width=0.8\textwidth]{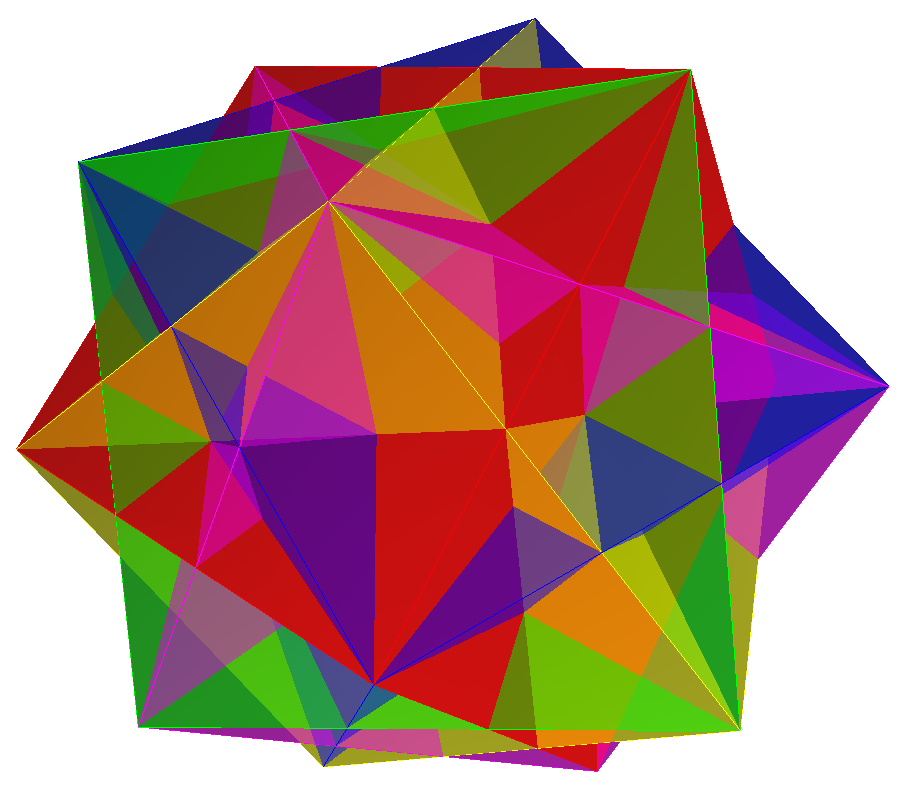}
\vspace{1.5cm}

\vfill
\vfill
{\Large Alexander Popa}

email: \href{mailto:alpopa@gmail.com}{alpopa@gmail.com}
\newpage
\textit{To my wife Raisa and my mother Tamara.}
\end{center}
\end{titlepage}
\tableofcontents
\newpage
\input{ch0-en.tex}
\input{ch1-en.tex}
\input{ch2-en.tex}
\input{ch3-en.tex}

\newpage
\thispagestyle{empty}
\mbox{}
\end{document}

%% file: ch0-en.tex
\chapter*{Introduction}
\addcontentsline{toc}{chapter}{Introduction}
\markboth{Introduction}{Introduction}
The first documented attempt to construct the geometry theory in an axiomatic way was made, as we know,
by Euclid (III cent BC) in his Elements. And while the word `geometry' literally means `earth measuring',
Euclidean geometry doesn't describe elliptic space, as more proper for measuring of our planet. New
axiomatic approach was revolutionary one, however the axiomatic has limitations. Euclid study what can
be constructed calculated or demonstrated starting with compass and straightedge. It was sufficient for
that time. However today, despite the fact Euclidean geometry is studied in the school, many people,
including geometriests, can't remember its axioms. Exception makes famous Euclid's V-th postulate,
which many of us remember in the form: ``At most one line can be drawn through any point not on a given
line parallel to the given line in a plane''. Euclid decided to formulate it so: ``If a line segment
intersects two straight lines forming two interior angles on the same side that sum to less than two
right angles, then the two lines, if extended indefinitely, meet on that side on which the angles sum
to less than two right angles''.

For modern geometry Euclid's axiomatic has several limitations:
\begin{itemize}
\item Euclid's axiomatic theory covered the only geometry system and only two--dimensional case. The
axiomatic of Euclidian geometry used today was developed by David Hilbert (1862 --- 1943), has 20 axioms
and covers two and three dimensions.
\item The four--dimensional case uses much more axioms. Development of axiomatic for spaces of further
dimensions is non--trivial.
\item Except hyperbolic geometry, construction of good axiomatic for other geometries is also non--trivial.
Usually, an axiomatic is constructed after the geometry is well studied with aim of some model (for example,
\cite{EWGL1912} describes the space--time axiomatic).
\item Undefined notions in geometry (point, line, between) differ very much from undefined notions
in other mathematic disciplines (number, function, space). Undefined notions of different geometries
differ from each other.
\item Mathematicians successful study Euclidean space of any dimension using analytic geometry and
forget Euclid's axiomatic.
\end{itemize}

Euclid's axiomatic played one important role. Its V-th postulate is so hard expressed and creates so
artificial feeling that urged mathematicians to create the hyperbolic geometry. Sad, when Nikolai
Lobachevsky (1792 --- 1856) and J\'anos Bolyai (1802 --- 1860) published their results, the new
geometry was slow in acceptance. Only after decades it was demonstrated that hyperbolic geometry is
interior geometry of surfaces with constant negative curvature. After next several years some models
of hyperbolic geometry were elaborated. Due to that fact the new geometry became accessible.

Author of a model, Felix Klein (1849 --- 1925) proposed ``Erlangen Program'' \cite{Erlangen} --- the
unified view over different geometries as complex of different transformation groups of space. The
invariants of these groups are figures of the geometries. In such way, Klein presented 9 two--dimensional
spaces. However, 6 of them he considered practic unaplicable \cite{Klein}. Till now speaking about
``non--euclidean geometry'', elliptic or hyperbolic geometry is primarily understood. Obviously, in order
to make all geometries to be taken seriously, an accessible model is required. One of such model for
two--dimensional case proposed \cite{Yaglom69} Isaak Moiseevich Yaglom (1921 --- 1988), using the notion
of generalized complex number. Among more recent results you can refer to \cite{Hach, Romakina, Artykbaev}.

In this work, supposed to your attention an uniform model of geometric spaces and based on it general
analytic geometry are described. Among its advantages there are its universality and linearity, hence
easyness to use. It isn't limited to specific dimension.

The first chapter describes different types of distance and angular measure and their models. Different
variants of axioms valid for different geometries are analyzed, as well as one variant of them, depending
on some parameter and universally valid. A analytic model depending on some parameters is constructed.
Lengths and angles are defined as parameters of corresponding motions.

In the second chapter you can find triangle equations valid for all geometries. The chapter describes
generalized orthogonal matrix as general form of motion matrix. A vector approach will be shown for
description of points, lines and planes, and for linear calculus of lengths and angles. At the end of
chapter, the reader will find a linear way to calculate volumes.

The third chapter has more philosophical character then practical one. Your attention will be set on
proper terminology and several well known spaces will be described in terms of constructed theory.

Uniform model of geometric spaces becomes the background of the GeomSpace project\footnote{\href
{http://sourceforge.net/projects/geomspace/}{http://sourceforge.net/projects/geomspace/}. The last version
of this book can also be downloaded form here.}.

%% file: ch1-en.tex
\chapter{Geometric Space Model Construction}
\section{Three Kinds of Plane Rotations. Rotation Characteristic}
Consider real plane $\mathbb{R}^2$. Consider three different transformations of $\mathbb{R}^2$:
rotation $\mathfrak{R}'(\phi)$, Galilean transformation $\mathfrak{R}''(\phi)$ and Lorintz
transformation $\mathfrak{R}'''(\phi)$ defined by matrices:
$$\mathfrak{R}'(\phi) = \begin{pmatrix}
\cos\phi&-\sin\phi\\
\sin\phi&\cos\phi
\end{pmatrix},$$
$$\mathfrak{R}''(\phi) = \begin{pmatrix}
1&0\\
\phi&1
\end{pmatrix},$$
and
$$\mathfrak{R}'''(\phi) = \begin{pmatrix}
\cosh\phi&\sinh\phi\\
\sinh\phi&\cosh\phi
\end{pmatrix},$$
where $\phi \in \mathbb{R}$.

\begin{figure}[h]
\begin{center}
\includegraphics{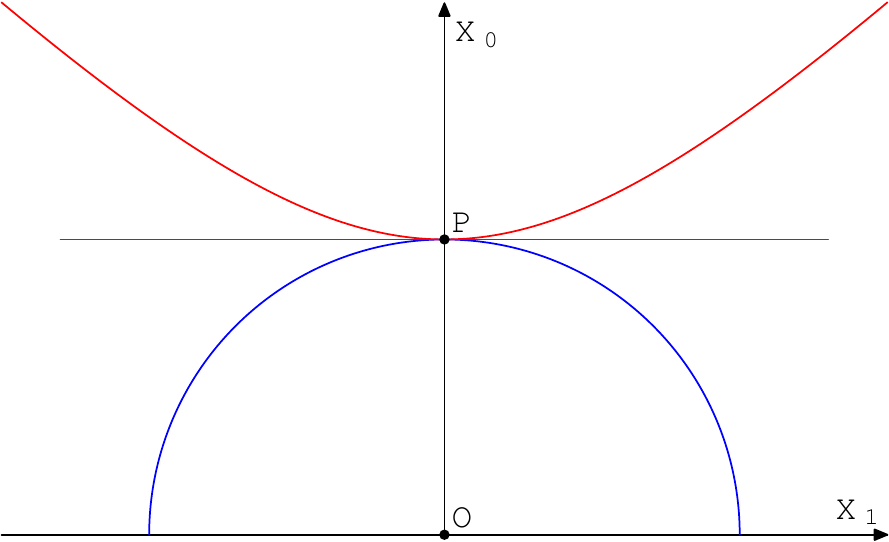}
\caption{Trajectory of point $P$ on transformations $\mathfrak{R}'$, $\mathfrak{R}''$ and $\mathfrak{R}'''$.}
\label{circles}
\end{center}
\end{figure}

Transformations $\mathfrak{R}'(\phi)$, $\mathfrak{R}''(\phi)$ and $\mathfrak{R}'''(\phi)$ have
several common properties. The determinant of all their matrices is 1, all them have the only fixed
point --- origin $O = (0, 0)$, $\mathfrak{R}(0) = I$ --- unit matrix, $\mathfrak{R}(x)\mathfrak{R}
(y) = \mathfrak{R}(x+y) = \mathfrak{R}(y)\mathfrak{R}(x)$ and the trajectory of point $P = (1, 0)$
verifies equations:
$$x_0^2 + x_1^2 = 1, \text{ for } \mathfrak{R}',$$
$$x_0 = 1, \text{ for } \mathfrak{R}'',$$
$$x_0^2 - x_1^2 = 1, \text{ for } \mathfrak{R}'''.$$
More general, the trajectory equation can be written as (Figure \ref{circles}):
$$x_0^2 + k\,x_1^2 = 1, k = -1, 0, 1.$$

We will name $\mathfrak{R}'$ \textit{elliptic rotation}, $\mathfrak{R}''$ \textit{parabolic rotation},
$\mathfrak{R}'''$ \textit{hyperbolic rotation} and $\phi$ respective angle. We will name the
coefficient $k$ \textit{characteristic} of a rotation. $k = 1$ corresponds to elliptic, $k = 0$ to
parabolic and $k = -1$ to hyperbolic rotation.

\section{Functions $C(x)$, $S(x)$ and $T(x)$}
We can see, that the matrices $\mathfrak{R}'$, $\mathfrak{R}''$ and $\mathfrak{R}'''$ have elements
$r_{11} = r_{22}$ and $r_{12} = -k\,r_{21}$. We can write:
\begin{equation}\label{rotation}
\mathfrak{R}(\phi)=\begin{pmatrix}
C(\phi)&-kS(\phi)\\
S(\phi)&C(\phi)
\end{pmatrix},
\end{equation}
where
$$
C(x) = \begin{cases}
\cos x,& k = 1\\
1,& k = 0\\
\cosh x,& k = -1
\end{cases}
$$
and
$$
S(x) = \begin{cases}
\sin x,& k = 1\\
x,& k = 0\\
\sinh x,& k = -1
\end{cases}
$$

Finally, we can define formally the functions $C(x)$ and $S(x)$ as\footnote{Here and further we will
consider for simplicity that $k^0 = 1$ for $k = 0$ too. We will say $x$ divide $k^i$, $k=0$ if in
expression $x/k^i$ the exponent of $k$ in numerator is not less then $i$.}:
\begin{equation}\label{cx}
C(x) = C(x, k) = \sum_{n=0}^\infty(-k)^n\frac{x^{2n}}{(2n)!},
\end{equation}
\begin{equation}\label{sx}
S(x) = S(x, k) = \sum_{n=0}^\infty(-k)^n\frac{x^{2n+1}}{(2n+1)!}.
\end{equation}
We will introduce one more function:
\begin{equation}\label{tx}
T(x) = \frac{S(x)}{C(x)}.
\end{equation}

Note, that always has place the equality:

\begin{equation}\label{cseq}
C^2(x) + kS^2(x) = 1,\, \forall x \in \mathbb{R}.
\end{equation}

We will name transformation (\ref{rotation}) \textit{generalized rotation}. Note, that along the
angle it has one more parameter --- its characteristic.

\section{Representation of Translation as Rotation. Its Characteristics}
Generalized rotation has a fixed point. However, the translation usually doesn't have a fixed
point\footnote{A fixed point is present for example in translation on elliptic plane.}. We will define
the translation through generalized rotation using an extra-dimension, as it is done in projective
geometry.

For $n$-dimensional space consider vector space $\mathbb{R}^{n+1}$. Rotation matrices have two
different rows compared to unit matrix. When one of them is the first row, we will consider them
translations. When none of them is the first one, we will consider them rotations. Therefore, the first
coordinate (we will count it 0) will be additional.

Let $n=1$. We will name the vector $o = \{1, 0\}$ origin. If $a = \mathfrak{R}(\phi)o$, then the
angle $\phi$ between $o$ and $a$ we can name distance $oa$. Different values of $k$ correspond to
different translations: when $k = 1$ translations are elliptic, when $k = 0$ they are parabolic, and
when $k = -1$ they are hyperbolic.

\begin{figure}[h]
\begin{center}
\includegraphics{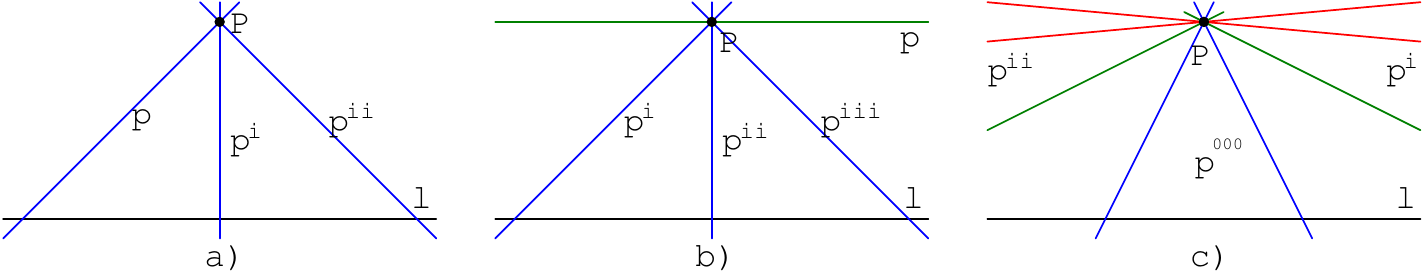}
\caption{Variants of Euclid's V postulate --- elliptic a), linear b) and hyperbolic c).}
\label{parallel}
\end{center}
\end{figure}

There is a kind of distance measure for each kind of translation: elliptic, parabolic and hyperbolic.
The difference between them can be seen in variants of V postulate of Euclid for elliptic, linear and
hyperbolic geometry (Figure \ref{parallel}).

Elliptic postulate (Figure \ref{parallel} a) is\footnote{For this case it is necessary to modify another
two postulates, namely that from any three points on a line exactly one lies between two others, and that
any line can be extended infinitely in any direction.}: For a given line $l$ and a point $P \notin l$,
exists no line $p \ni P\,|\, l\cap p = \varnothing$. It is identical to the following: For a given line
$l$ and a point $P \notin l$, all lines $p \ni P$ intersect $l$.

The linear postulate (Figure \ref{parallel} b) is: For a given line $l$ and a point $P \notin l$, exists
one line $p \ni P\,|\, l\cap p = \varnothing$.

The hyperbolic postulate (Figure \ref{parallel} c) is: For a given line $l$ and a point $P \notin l$,
exist at least two lines $p', p''\ni P\,|\, l\cap p' = \varnothing, l\cap p'' = \varnothing$.

Generally, V postulate of Euclid can be formulated as: For a given line $l$ and a point $P \notin l$,
exist $0^{k_1}$ lines $p \ni P\,|\, l\cap p = \varnothing$. It should be mentioned that $0^{k_1}$ is a
symbol, not a number used in calculus. Its value equals to 0 for $k_1 = 1$, 1 for $k_1 = 0$ and $\infty$
for $k_1 = -1$.

\section{Kinds of Space Rotations. Bundles of Unconnectable Points}
\begin{figure}[h]
\begin{center}
\includegraphics{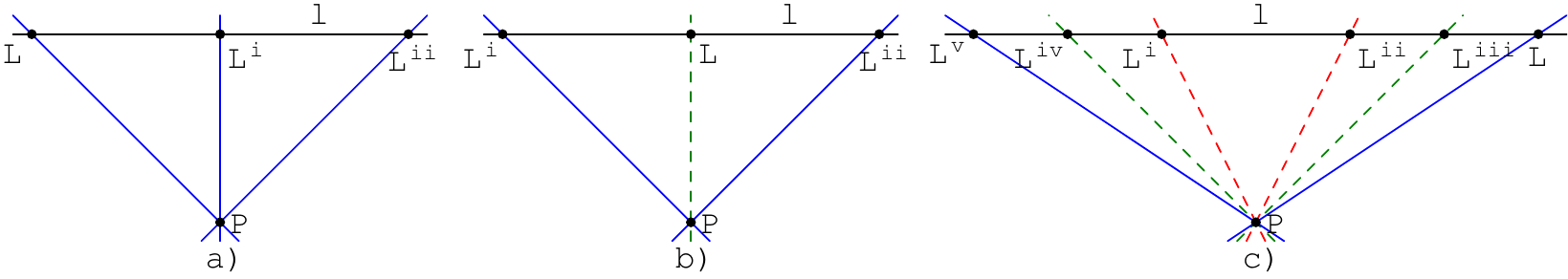}
\caption{Different variants of points unconnectability axiom --- elliptic a), linear b) and hyperbolic c).}
\label{connectable}
\end{center}
\end{figure}

It's easy to see that classic rotations in Euclidean geometry, as well as in the elliptic (Riemannian)
geometry and the hyperbolic (Bolyai--Lobachevsky) geometry has the characteristic $k = 1$. We can extend
the notion of space rotation to generalized space rotation with some characteristic. The best way to
illustrate difference between them is to formulate angular equivalent of V Postulate of Euclid --- axiom
of points connectability (Figure \ref{connectable}). In order to do this we will change the following
phrases between them:
\begin{eqnarray*}
\text{line } l &\longleftrightarrow& \text{ point } L,\\
P \in l &\longleftrightarrow& p \ni L,\\
P \notin l &\longleftrightarrow& p \not\ni L,\\
AB = \phi &\longleftrightarrow& \angle ab = \phi,\\
a \cap b = C &\longleftrightarrow& c = AB,\\
a \cap b = \varnothing &\longleftrightarrow& A \text{ is unconnectable with } B.
\end{eqnarray*}

The last statement is unusual for the above three geometries\footnote{It conflicts with axiom which
states that through any two points goes a line. This axiom should be changed by one of the following in
order to consider the geometries with non-elliptic rotations.}. It makes sense in geometries with angular
characteristic 0 or $-1$. The unconnectable property of points is similar to parallel property of lines.

The angle equivalent of V Postulate for elliptic characteristic (Figure \ref{connectable} a) is: On a
line $l \not\ni P$ exist no points $L$ unconnectable with $P$.

For parabolic characteristic (Figure \ref{connectable} b) it is: On a line $l \not\ni P$ exists the only
point $L$ unconnectable with $P$.

For the hyperbolic characteristic (Figure \ref{connectable}) it is: On a line $l \not\ni P$ exist at
least two points $L'$ and $L''$ unconnectable with $P$.

Generally this axiom can be formulated as: On a line $l \not\ni P$ exist $0^{k_2}$ points $L$
unconnectable with $P$. As in case of parallel lines, symbol $0^{k_2}$ isn't used in calculus.

Similar to bundles of lines --- intersected, parallel or divergent we can speak about \textit
{bundles of points}. More exactly, let $X, Y \in \mathbb{R}^{n+1}$. All linear combinations
$Z = \alpha X + \beta Y, \alpha, \beta \in \mathbb{R}$ form a set we will name bundle of points.
As we will see, this set has one constraint. Therefore, it has one free parameter. As every two lines
define a bundle of lines, every two points ($X$ and $Y$) define bundle of points. If $X$ is
connectable with $Y$ this bundle is a line (similar to intersection point of bundle of
intersected lines). Lines has blue color on figure \ref{connectable}. If $X$ and $Y$ are unconnectable,
this bundle of points isn't a line (similar to bundle of parallel or divergent lines). Bundles of
unconnectable points are green and red on figure \ref{connectable}.

\begin{figure}[h]
\begin{center}
\includegraphics{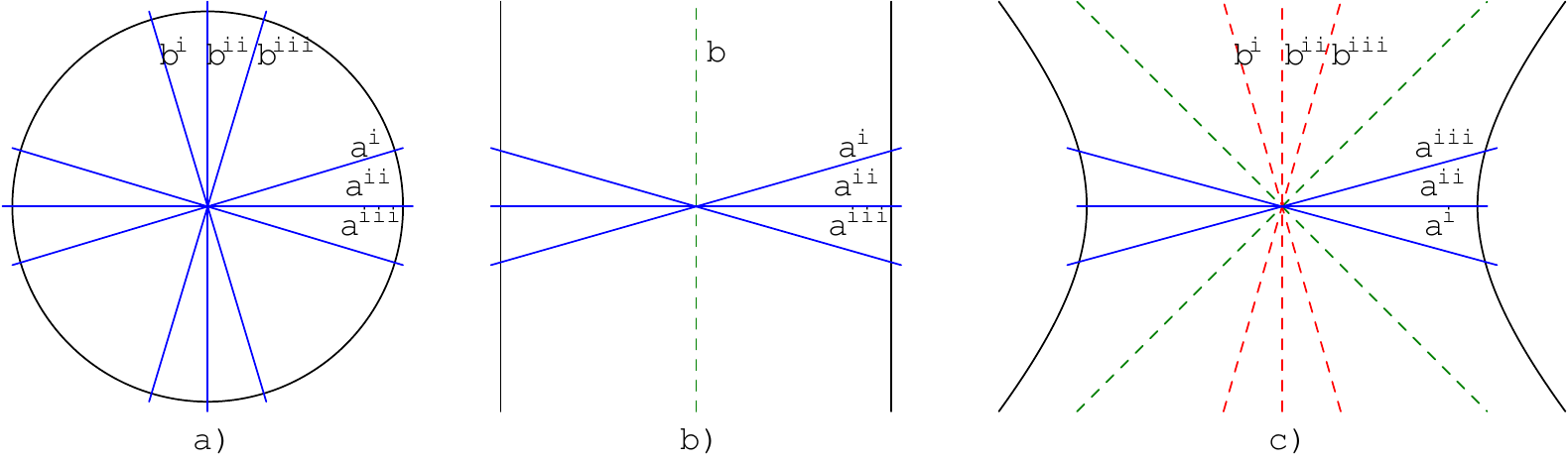}
\caption{Mutual position of different bundles and circles a) elliptic angular characteristic, b) linear
angular characteristic and c) hyperbolic angular characteristic.}
\label{bundles}
\end{center}
\end{figure}

For any angle characteristic there are infinity of bundles of connectable points. For angle
characteristic 1 all point bundles are lines. For angle characteristic 0 for any point there is the
only bundle of unconnectable points (green). For angle characteristic $-1$ there are infinity bundles of
unconnectable points (red). In thes case the bundles of connectable points and the bundles of
unconnectable points for some point form two categories of bundles. The limit (marginal) bundles
of unconnectable points (green) can be viewed as the third category (similar to differencee between
parallel and divergent lines). There are exactly two limit bundles. Note that bundles of connectable
points intersect all circles with centre in the centre of bundle, all bundles of unconnectable points
don't intersect these circles and limit bundles are asymptotic to circles (Figure \ref{bundles}).

Emphasize that the angle between two lines and the angle between two two--dimensional planes are the
different measures. The angle between two threedimensional planes is different from them both and so
on. Thus, the angle between lines can have the different characteristic then the angle between
two--dimensional planes and so on.

\section{Main Space Rotations}
Consider $\mathbb{R}^{n+1}$ and $k_1, k_2,... k_n \in \{-1, 0, 1\}$. We will note $C_i(x) = C(x, k_i)$,
$S_i(x) = S(x, k_i)$ and $T_i(x) = S_i(x)/C_i(x)$. Let
$$\mathfrak{R}_1(\phi)=\begin{pmatrix}
C_1(\phi)&-k_1S_1(\phi)&0&\ldots&0\\
S_1(\phi)&C_1(\phi)&0&\ldots&0\\
0&0&1&\ldots&0\\
\vdots&\vdots&\vdots&\ddots&\vdots\\
0&0&0&\ldots&1
\end{pmatrix},$$
$$\mathfrak{R}_2(\phi)=\begin{pmatrix}
1&0&0&\ldots&0\\
0&C_2(\phi)&-k_2S_2(\phi)&\ldots&0\\
0&S_2(\phi)&C_2(\phi)&\ldots&0\\
\vdots&\vdots&\vdots&\ddots&\vdots\\
0&0&0&\ldots&1
\end{pmatrix},$$
$$\vdots$$
$$\mathfrak{R}_n(\phi)=\begin{pmatrix}
1&0&\ldots&0&0\\
0&1&\ldots&0&0\\
\vdots&\vdots&\ddots&\vdots&\vdots\\
0&0&\ldots&C_n(\phi)&-k_nS_n(\phi)\\
0&0&\ldots&S_n(\phi)&C_n(\phi)
\end{pmatrix}.$$

We will name $\mathfrak{R}_1,... \mathfrak{R}_n$ \textit{main space rotations}.

\section{Vector Product. Invariant Quadric Form}
Let
\begin{equation}\label{K}
K_m = \prod_{i=1}^mk_i,\,\forall m = \overline{0,n}
\end{equation}
We can see that $K_m \in \{-1,0,1\},\, \forall m=\overline{0,n}$ as well as $k_m$. Let define vector
product $\odot$ as
\begin{equation}\label{vdot}
x\odot y = \sum_{i=0}^nK_ix_iy_i
\end{equation}

For some vectors $x = \{x_0, x_1,... x_n\}$ and $y = \{y_0, y_1,... y_n\}$, $x' = \mathfrak{R}_m(\phi)x
= \{x_0,... x_{m-2}, x_{m-1}$ $C_m(\phi) - k_m x_m S_m(\phi), x_{m-1} S_m(\phi)+x_m C_m(\phi), x_{m+1},...
x_n\}$ and $y' = \mathfrak{R}_m(\phi)y$. We can see that
\begin{eqnarray*}
x'\odot y'&=&\sum_{i=0}^nK_ix'_iy'_i\\
&=& \sum_{i=0}^{m-2}K_ix_iy_i\\
&+& ((x_{m-1}C_m(\phi)-k_mx_mS_m(\phi))(y_{m-1}C_m(\phi)-k_my_mS_m(\phi))\\
&+& k_m(x_{m-1}S_m(\phi)+x_mC_m(\phi))(y_{m-1}S_m(\phi)+y_mC_m(\phi)))K_{m-1}\\
&+& \sum_{i=m+1}^nK_ix_iy_i\\
&=& \sum_{i=0}^{m-2}K_ix_iy_i\\
&+& ((x_{m-1}y_{m-1}C_m^2(\phi) - k_m(x_{m-1}y_m + x_my_{m-1})S_m(\phi)C_m(\phi)
+ k_m^2x_my_mS_m^2(\phi)\\
&+& k_mx_{m-1}y_{m-1}S_m^2(\phi) + k_m(x_{m-1}y_m + x_my_{m-1})S_m(\phi)C_m(\phi)
+ k_mx_my_mC_m^2(\phi)))K_{m-1}\\
&+& \sum_{i=m+1}^nK_ix_iy_i\\
&=& \sum_{i=0}^{m-2}K_ix_iy_i\\
&+& (x_{m-1}y_{m-1}(C_m^2(\phi) + k_mS_m^2(\phi)) + k_mx_my_m(C_m^2(\phi) + k_mS_m^2(\phi)))K_{m-1}\\
&+& \sum_{i=m+1}^nx_iy_iK_i\\
&=& \sum_{i=0}^nK_ix_iy_i = x\odot y
\end{eqnarray*}
This is true for all $m = \overline{1,n}$. So the quadric form $x\odot y$ is invariant in respect to
main rotations of $\mathfrak{R}_m$.

\section{Space Definition by its Specification}
\begin{figure}[h]
\begin{center}
\includegraphics{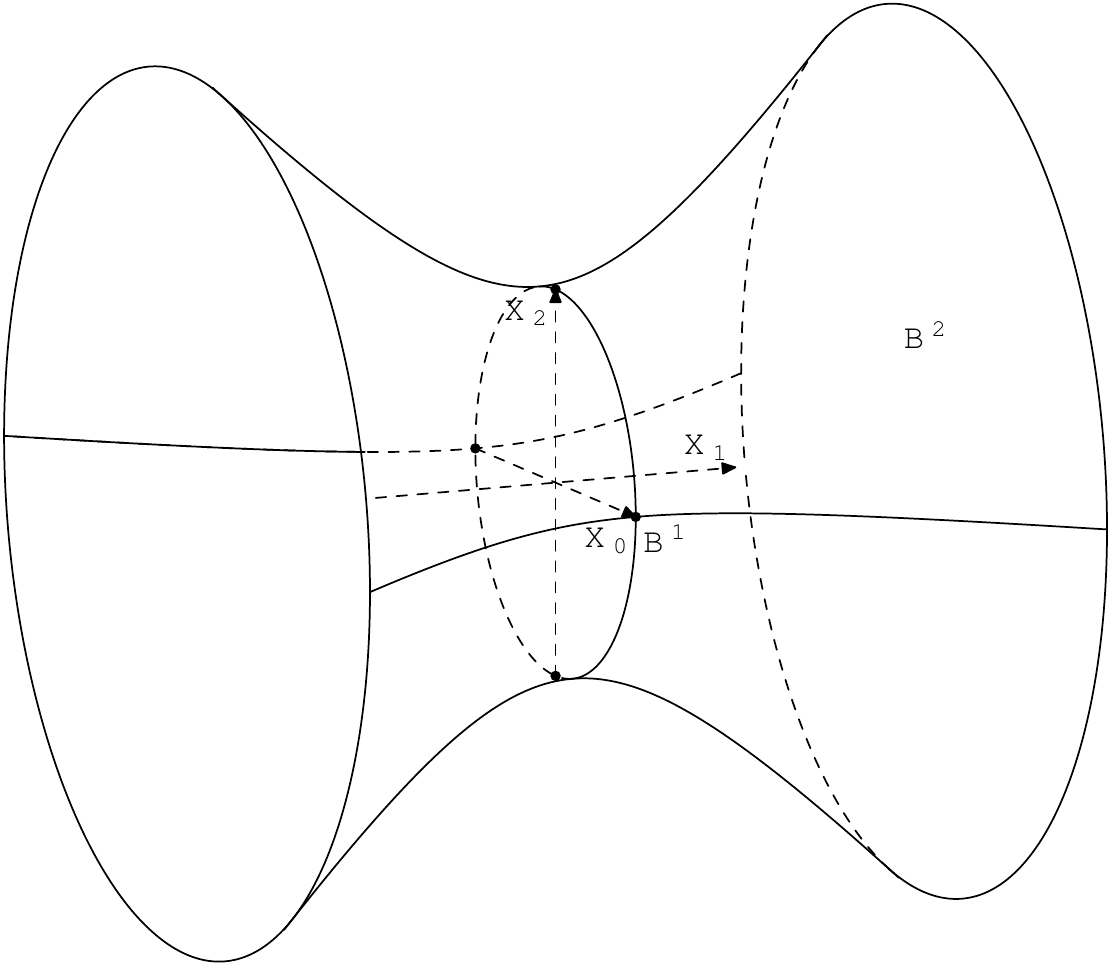}
\caption{Sphere of space with specification $\{-1, -1\}$.}
\label{sphere}
\end{center}
\end{figure}

Consider $\mathbb{RP}^n$ projective space and $k_i \in \{-1,0,1\},\, \forall i = \overline{1,n}$. We
can now introduce a geometric space `unit sphere' $\mathbb{B}^n = \{x \in \mathbb{RP}^n\, |\, x\odot x =
1\}$ (Figure \ref{sphere}). As all main rotations preserves the quadric form defined by product $\odot$,
they also preserves $\mathbb{B}^n$. We will name $k_i, i = \overline{1,n}$ \textit{space specification}.
We will name `point' $X \in \mathbb{B}^n$ the corresponding vector $x \in \mathbb{RP}^n$ and will use
homogeneous coordinates normalized in order to $x\odot x = 1$.

We will name `origin' of $\mathbb{B}^n$ the point $O = \left[1:0:...:0\right] \in \mathbb{B}^n$. It
isn't origin of $\mathbb{R}^{n+1}$, $(0, 0,... 0)\notin \mathbb{B}^n$ and we will refer to $O$ as
origin if isn't specified otherwise.

It's easy to see that for any $k_1, k_2,... k_n$, $O = \mathbb{B}^0 \subset \mathbb{B}^1 \subset ...
\subset \mathbb{B}^n$.

We will define motions of $\mathbb{B}^n$ all transformations that result on finite product of main
rotations.

We will define `lines' all images of $\mathbb{B}^1$ on any motion of $\mathbb{B}^n$. Similarly,
we define `$m$-dimensional' planes all images of $\mathbb{B}^m$ on any motions of $\mathbb{B}^n$
for any $m \in \overline{0,n-1}$.

For each characteristic parameter $k_i$ we can introduce a scale parameter $r_i \in \mathbb{R}_+$,
$i = \overline{1,n}$. The $k_1/r_1^2$ is exactly the gaussian curvature of space. Others have no
representation since finite angle measure doesn't require scaling. In this case the radian measure is
native. An example of angle scale is degree measure which has scale $180/\pi$. However when the angle
is not bounded a scale introduction has sense. All scales can be easy embedded in functions $C_i(x)$,
$S_i(x)$ and $T_i(x)$ by using instead $C_i\left(\frac{x}{r_i}\right)$, $S_i\left(\frac{x}{r_i}
\right)$ and $T_i\left(\frac{x}{r_i}\right)$ respectively, $i = \overline{1,n}$.

\section{Definition of Measure Using Motions}
A traditional way of definition the measures and motions is to provide a way to calculate the
distances as is and then to define motions in such way that all maps $M:\mathbb{R}^n\to\mathbb
{R}^n$ preserve the distance. We go another way. We provide motions as is and then search for a way
to define measures in such way that motions preserve them.

We will say point $A \in \mathbb{B}^1\subset\mathbb{B}^n$ has the distance $\phi$ from origin $O$ if
$A = \mathfrak{R}_1(\phi)O$. Having $O = \left[1:0:...:0\right]$, $A = \left[C_1(\phi):S_1(\phi):0:
...:0\right]$, $O\odot A = C_1(\phi)$. We will say one--dimensional (planar) angle between $\mathbb{B}^1$
and some one--dimensional line $\mathbb{B}'^1 \subset \mathbb{B}^2$ equals $\phi$ if $\mathbb{B}'^1 =
\mathfrak{R}_2(\phi)\mathbb{B}^1$. Similarly, we will define the $m$-dimensional angle $\phi$ between
$\mathbb{B}^m$ and $m$-dimensional plane $\mathbb{B}'^m \subset \mathbb{B}^{m+1}$ if $\mathbb{B}'^m =
\mathfrak{R}_{m+1}(\phi)\mathbb{B}^m,\, \forall m = \overline{0, m-1}$. Note, that $n$-dimensional
angle between any planes is 0 since all them are subset of $B^n$.

Let $X, Y \in \mathbb{B}^n$. If there exists a motion that maps $\mathbb{B}^1$ to $XY$ we will
name points $X$ and $Y$ connectable and distance $XY$ measurable. If not, we will name points $X$ and
$Y$ unconnectable (just as lines can be parallel) and strictly speaking the distance $XY$ doesn't
exists\footnote{In this case there exist a measure $XY$, but it may have different characteristic then
distance. We will name this measure also distance, keeping in mind that it is generalized distance.}.

We can find a motion $\mathfrak{M}$ of space $\mathbb{B}^n$ that maps origin $O$ to $X$ and some
point $A \in\mathbb{B}^1\subset\mathbb{B}^n$ to $Y$. As motion $\mathfrak{M}$ preserves the quadric
form $\odot$, we can see that $X\odot Y = O\odot A$. We can define the distance $\phi$ between $X$ and
$Y$ as
\begin{equation}\label{c1meas}
C_1(\phi) = X\odot Y.
\end{equation}

It's easy to see that all motions preserve the distance. In case of elliptic, Euclidian and
hyperbolic space it is sufficient, because all other measures can be calculated from distances. However,
in some spaces angles can be scaled in a manner distances are scaled in Euclidean space. So we should
find the way to measure all the measures in general case.

%% file: ch2-en.tex
\chapter{Measure Calculus}
\section{General Triangle Equations}
\begin{figure}[h]
\begin{center}
\includegraphics{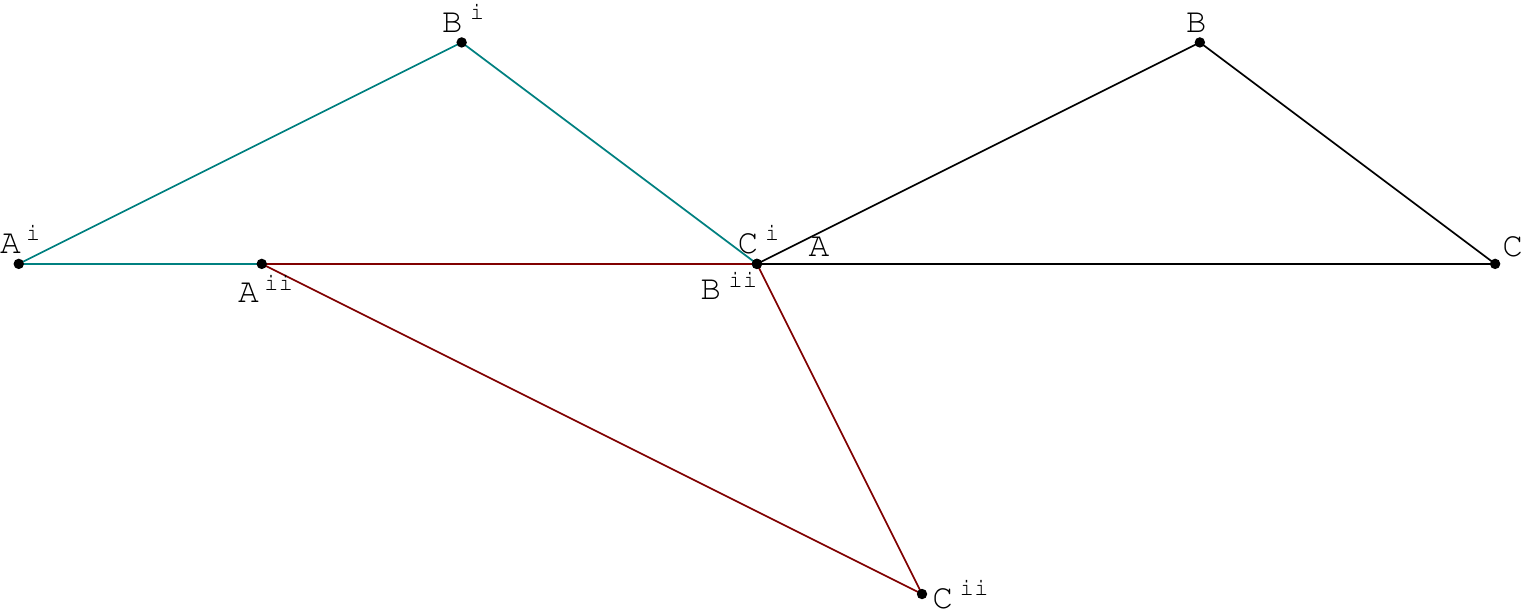}
\caption{General triangle equations deduction.}
\label{general.triangle}
\end{center}
\end{figure}

Consider triangle $ABC \in \mathbb{B}^2$ with the edges $a$, $b$, $c$, interior angles $\alpha$,
$\gamma$ and exterior angle $\beta'$ (Figure \ref{general.triangle}). Let $A = \left[1:0:0\right]$
the origin, $C = \mathfrak{R}_1(b)A = \left[C_1(b):S_1(b):0\right]$ and $B = \mathfrak{R}_2(\alpha)
\mathfrak{R}_1(c)A = \left[C_1(c):S_1(c)C_2(\alpha):S_1(c)S_2(\alpha)\right]$. Note, that the
interior angle $\beta$ does not exist in case of $k_2 = 0$ or $k_2 = -1$. The exterior angle
$\beta'$ always exists.

Now, let $A'B'C' = \mathfrak{R}(-b)(ABC)$ (Figure \ref{general.triangle}, cyan). The point we are
interested in is $B' = \left[C_1(b) C_1(c) + k_1 S_1(b) S_1(c) C_2(\alpha): -S_1(b) C_1(c) +
C_1(b) S_1(c) C_2(\alpha): S_1(c) S_2(\alpha)\right]$. At the other hand, now  $B' = \left[C_1(a):
-S_1(a) C_2(\gamma): S_1(a) S_2(\gamma)\right]$. It means:
\begin{eqnarray*}
C_1(a) &=& C_1(b)C_1(c) + k_1S_1(b)S_1(c)C_2(\alpha),\\
-S_1(a)C_2(\gamma) &=& -S_1(b)C_1(c) + C_1(b)S_1(c)C_2(\alpha),\\
S_1(a)S_2(\gamma) &=& S_1(c)S_2(\alpha).
\end{eqnarray*}
The first equation is the form of the Cosine I law. Similarly we have
$$C_1(c) = C_1(a)C_1(b) + k_1S_1(a)S_1(b)C_2(\gamma).$$
The third equation is equivalent to
$$\frac{S_1(a)}{S_2(\alpha)} = \frac{S_1(c)}{S_2(\gamma)},$$
which is the form of the Sine law.

Let now $A''B''C'' = \mathfrak{R}_1(-c)\mathfrak{R}_2(-\alpha)(ABC)$ (Figure \ref{general.triangle},
brown). Now we are interested in vertex $C'' = \left[C_1(b) C_1(c) + k_1 S_1(b) S_1(c) C_2(\alpha):
-C_1(b) S_1(c) + S_1(b) C_1(c) C_2(\alpha): -S_1(b) S_2(\alpha)\right]$. At the other hand, $C''=
\left[C_1(a): S_1(a) C_2(\beta'): -S_1(a) S_2(\beta')\right]$. From  here we have:
\begin{eqnarray*}
C_1(a) &=& C_1(b)C_1(c) + k_1S_1(b)S_1(c)C_2(\alpha),\\
S_1(a)C_2(\beta') &=& -C_1(b)S_1(c) + S_1(b)C_1(c)C_2(\alpha),\\
-S_1(a)S_2(\beta') &=& -S_1(b)S_2(\alpha).
\end{eqnarray*}
The first one is the Cosine I law, the third one is equivalent to:
\begin{equation}\label{sinlaw}
\frac{S_1(a)}{S_2(\alpha)} = \frac{S_1(b)}{S_2(\beta')} = \frac{S_1(c)}{S_2(\gamma)}
\end{equation}
which is the Sine law. Note that in case $k_2 = 1$ we have $\beta = \pi - \beta'$,
$S_2(\beta) = S_2(\beta')$. Let calculate the value of $C_2(\alpha)$ from the first equation and
put it to the second one:
\begin{eqnarray*}
S_1(a)C_2(\beta') &=& -C_1(b)S_1(c) + S_1(b)C_1(c)\frac{C_1(a) - C_1(b)C_1(c)}{k_1S_1(b)S_1(c)}\\
&=& -C_1(b)S_1(c) + C_1(c)\frac{C_1(a) - C_1(b)C_1(c)}{k_1S_1(c)},\\
k_1S_1(a)S_1(c)C_2(\beta') &=& -k_1S_1(c)^2C_1(b) + C_1(a)C_1(c) - C_1(b)C_1(c)^2\\
&=& C_1(a)C_1(c) - C_1(b)(C_1(c)^2 + k_1S_1(c)^2)\\
&=& C_1(a)C_1(c) - C_1(b),\\
C_1(b) &=& C_1(a)C_1(c) - k_1S_1(a)S_1(c)C_2(\beta').
\end{eqnarray*}
Note the `$-$' sign in the right part of the equation. It is so because the $\beta'$ angle is external.
For the case $k_2 = 1$, the internal angle $\beta = \pi - \beta'$, $C_2(\beta) = -C_2(\beta')$.

What about the Cosine II law? We will use these two equations:
\begin{eqnarray*}
-S_1(a)C_2(\gamma) &=& -S_1(b)C_1(c) + C_1(b)S_1(c)C_2(\alpha),\\
S_1(a)C_2(\beta') &=& -C_1(b)S_1(c) + S_1(b)C_1(c)C_2(\alpha)
\end{eqnarray*}
First, replace $S_1(b)$ with $S_1(a)S_2(\beta')/S_2(\alpha)$ and $S_1(c)$ with
$S_1(a)S_2(\gamma)/S_2(\alpha)$:
\begin{eqnarray*}
-S_1(a)C_2(\gamma) &=& -S_1(a)\frac{S_2(\beta)}{S_2(\alpha)}C_1(c)
+ C_1(b)S_1(a)\frac{S_2(\gamma)}{S_2(\alpha)}C_2(\alpha),\\
-S_2(\alpha)C_2(\gamma) &=& -C_1(c)S_2(\beta') + C_1(b)S_2(\gamma)C_2(\alpha),\\
S_2(\beta')C_1(c) &=& S_2(\alpha)C_2(\gamma) + C_2(\alpha)S_2(\gamma)C_1(b),
\end{eqnarray*}
and
\begin{eqnarray*}
S_1(a)C_2(\beta') &=& -C_1(b)S_1(a)\frac{S_2(\gamma)}{S_2(\alpha)}
+ S_1(a)\frac{S_2(\beta')}{S_2(\alpha)}C_1(c)C_2(\alpha),\\
S_2(\alpha)C_2(\beta') &=& -C_1(b)S_2(\gamma) + C_1(c)S_2(\beta')C_2(\alpha),\\
S_2(\gamma)C_1(b) &=& -S_2(\alpha)C_2(\beta') + C_2(\alpha)S_2(\beta')C_1(c).
\end{eqnarray*}
Now from the first equation let calculate $C_1(c)$ and put it in the second one:
\begin{eqnarray*}
S_2(\gamma)C_1(b) &=& -S_2(\alpha)C_2(\beta') + C_2(\alpha)S_2(\beta')
\frac{S_2(\alpha)C_2(\gamma) + C_2(\alpha)S_2(\gamma)C_1(b)}{S_2(\beta')}\\
&=& -S_2(\alpha)C_2(\beta') + C_2(\alpha)S_2(\alpha)C_2(\gamma) + C_2(\alpha)^2S_2(\gamma)C_1(b),\\
S_2(\gamma)C_1(b)(1 - C_2(\alpha)^2) &=& S_2(\alpha)(C_2(\alpha)C_2(\gamma) - C_2(\beta')),\\
k_2S_2(\gamma)C_1(b)S_2(\alpha)^2 &=& S_2(\alpha)(C_2(\alpha)C_2(\gamma) - C_2(\beta)),\\
k_2S_2(\alpha)S_2(\gamma)C_1(b) &=& C_2(\alpha)C_2(\gamma) - C_2(\beta'),\\
C_2(\beta') &=& C_2(\alpha)C_2(\gamma) - k_2S_2(\alpha)S_2(\gamma)C_1(b).
\end{eqnarray*}
When $k_2 = 1$ we have
\begin{eqnarray*}
-\cos\beta &=& \cos\alpha\cos\gamma - k_2\sin\alpha\sin\gamma C_1(b),\\
\cos\beta &=& -\cos\alpha\cos\gamma + k_2\sin\alpha\sin\gamma C_1(b).
\end{eqnarray*}
Similarly, calculating $C_1(b)$ form the second equation and putting it in the first one, obtain:
\begin{eqnarray*}
S_2(\beta')C_1(c) &=& S_2(\alpha)C_2(\gamma) + C_2(\alpha)S_2(\gamma)
\frac{C_2(\alpha)S_2(\beta')C_1(c) - S_2(\alpha)C_2(\beta')}{S_2(\gamma)}\\
&=& S_2(\alpha)C_2(\gamma) + C_2(\alpha)^2S_2(\beta')C_1(c) - C_2(\alpha)S_2(\alpha)C_2(\beta'),\\
S_2(\beta')C_1(c)(1 - C_2(\alpha)^2) &=& S_2(\alpha)(C_2(\gamma) - C_2(\alpha)C_2(\beta')),\\
k_2S_2(\beta')C_1(c)S_2(\alpha)^2 &=& S_2(\alpha)(C_2(\gamma) - C_2(\alpha)C_2(\beta')),\\
k_2S_2(\alpha)S_2(\beta')C_1(c) &=& C_2(\gamma) - C_2(\alpha)C_2(\beta'),\\
C_2(\gamma) &=& C_2(\alpha)C_2(\beta') + k_2S_2(\alpha)S_2(\beta')C_1(c).
\end{eqnarray*}
When $k_2 = 1$ we have as above
$$\cos\gamma = -\cos\alpha\cos\beta + k_2\sin\alpha\sin\beta C_1(c).$$
Similarly, we have
$$C_2(\alpha) = C_2(\beta')C_2(\gamma) + k_2S_2(\beta')S_2(\gamma)C_1(a).$$

We will find the form of the Cosine I and II law that does not contain $C_1$ or $C_2$ functions in the
left part. However, it contains these functions in the right part. It makes sense since in the case
$k_1 \ne 0$ (for the Cosine I law) and $k_2 \ne 0$ (for the Cosine II law) when we can calculate their
respective $C^{-1}$ functions, but when $k_1 = 0$, the space admit distance scaling and the angle values
does not determine the distances (Cosine II law is a equality which doesn't contain $C_1$ function),
while when $k_2 = 0$, the space admit the angular scaling and distances does not determine angles (Cosine
I law is a equality which doesn't contain $C_2$ function).

Note also that we can deduce one form of Cosine I and one form of Cosine II law if we introduce a (may be
virtual) angle $\beta$ so as:
\begin{eqnarray*}
S_2(\beta) &=& S_2(\beta'),\\
C_2(\beta) &=& -C_2(\beta'),\\
T_2(\beta) &=& -T_2(\beta').
\end{eqnarray*}
Then both Cosine I and II law have identical form. Now let calculate
\begin{eqnarray*}
k_1S_1^2(a) &=& 1 - C_1^2(a)\\
&=& (C_1^2(b) + k_1S_1^2(b))(C_1^2(c) + k_1S_1^2(c))\\
&-& (C_1(b)C_1(c) + k_1S_1(b)S_1(c)C_2(\alpha))^2\\
&=& C_1^2(b)C_1^2(c) + k_1C_1^2(b)S_1^2(c) + k_1S_1^2(b)C_1^2(c) + k_1^2S_1^2(b)S_1^2(c)\\
&-& C_1^1(b)C_1^2(c) - 2k_1C_1(b)C_1(c)S_1(b)S_1(c)C_2(\alpha) - k_1^2S_1^2(b)S_1^2(c)C_2^2(\alpha)\\
&=& k_1(C_1^2(b)S_1^2(c) + S_1^2(b)C_1^2(c) - 2C_1(b)C_1(c)S_1(b)S_1(c)C_2(\alpha))\\
&+& k_1^2S_1^2(b)S_1^2(c)(1 - C_2^2(\alpha))\\
&=& k_1(C_1^2(b)S_1^2(c) + S_1^2(b)C_1^2(c) - 2C_1(b)C_1(c)S_1(b)S_1(c)C_2(\alpha))\\
&+& k_1^2k_2S_1^2(b)S_1^2(c)S_2^2(\alpha),\\
S_1^2(a) &=& C_1^2(b)S_1^2(c) + S_1^2(b)C_1^2(c) - 2C_1(b)C_1(c)S_1(b)S_1(c)C_2(\alpha)\\
&+& k_1k_2S_1^2(b)S_1^2(c)S_2^2(\alpha),
\end{eqnarray*}
or, having:
\begin{eqnarray}
C_1(a) &=& C_1(b)C_1(c)(1 + k_1T_1(b)T_1(c)C_2(\alpha)),\nonumber\\
T_1^2(a) &=& \frac{T_1^2(b) + T_1^2(c) - 2T_1(b)T_1(c)C_2(\alpha) + k_1k_2T_1^2(b)T_1^2(c)S_1^2(\alpha)}
{(1 + k_1T_1(b)T_1(c)C_2(\alpha))^2}.\label{cos1a}
\end{eqnarray}
Similarly,
\begin{eqnarray}
S_1^2(b) &=& C_1^2(a)S_1^2(c) + S_1^2(a)C_1^2(c) + 2C_1(a)C_1(c)S_1(a)S_1(c)C_2(\beta')\nonumber\\
&+& k_1k_2S_1^2(a)S_1^2(c)S_2^2(\beta'),\nonumber\\
T_1^2(b) &=& \frac{T_1^2(a) + T_1^2(c) + 2T_1(a)T_1(c)C_2(\beta') + k_1k_2T_1^2(a)T_1^2(c)S_1^2(\beta')}
{(1 - k_1T_1(a)T_1(c)C_2(\beta'))^2},\label{cos1b}
\end{eqnarray}
and
\begin{eqnarray}
S_1^2(c) &=& C_1^2(a)S_1^2(b) + S_1^2(a)C_1^2(b) - 2C_1(a)C_1(b)S_1(a)S_1(b)C_2(\gamma)\nonumber\\
&+& k_1k_2S_1^2(a)S_1^2(b)S_2^2(\gamma),\nonumber\\
T_1^2(c) &=& \frac{T_1^2(a) + T_1^2(b) - 2T_1(a)T_1(b)C_2(\gamma) + k_1k_2T_1^2(a)T_1^2(b)S_1^2(\gamma)}
{(1 + k_1T_1(a)T_1(b)C_2(\gamma))^2}.\label{cos1c}
\end{eqnarray}

Now, let calculate
\begin{eqnarray*}
k_2S_2^2(\alpha) &=& 1 - C_2^2(\alpha)\\
&=& (C_2^2(\beta') + k_2S_2^2(\beta'))(C_2^2(\gamma) + k_2S_2^2(\gamma))\\
&-& (C_2(\beta')C_2(\gamma) + k_2S_2(\beta')S_2(\gamma)C_1(a))^2\\
&=& C_2^2(\beta')C_2^2(\gamma) + k_2C_2^2(\beta')S_2^2(\gamma)
+ k_2S_2^2(\beta')C_2^2(\gamma) + k_2^2S_2^2(\beta')S_2^2(\gamma)\\
&-& C_2^2(\beta')C_2^2(\gamma) - 2k_2C_2(\beta')S_2(\beta')C_2(\gamma)S_2(\gamma)C_1(a)
- k_2^2S_2^2(\beta')S_2^2(\gamma)C_1^2(a)\\
&=& k_2(C_2^2(\beta')S_2^2(\gamma) + S_2^2(\beta')C_2^2(\gamma)
- 2C_2(\beta')S_2(\beta')C_2(\gamma)S_2(\gamma)C_1(a))\\
&+& k_2^2S_2^2(\beta')S_2^2(\gamma)(1 - C_1^2(a))\\
&=& k_2(C_2^2(\beta')S_2^2(\gamma) + S_2^2(\beta')C_2^2(\gamma)
- 2C_2(\beta')S_2(\beta')C_2(\gamma)S_2(\gamma)C_1(a))\\
&+& k_1k_2^2S_2^2(\beta')S_2^2(\gamma)S_1^2(a),\\
S_2^2(\alpha) &=& C_2^2(\beta')S_2^2(\gamma) + S_2^2(\beta')C_2^2(\gamma)
- 2C_2(\beta')S_2(\beta')C_2(\gamma)S_2(\gamma)C_1(a)\\
&+& k_1k_2S_2^2(\beta')S_2^2(\gamma)S_1^2(a),
\end{eqnarray*}
or
\begin{equation}\label{cos2a}
T_2^2(\alpha) = \frac{T_2^2(\beta') + T_2^2(\gamma) - 2T_2(\beta')T_2(\gamma)C_1(a)
+ k_1k_2T_2^2(\beta')T_2^2(\gamma)S_1^2(a)}{(1 + k_2T_2(\beta')T_2(\gamma)C_1(a))^2}.
\end{equation}
Similarly,
\begin{eqnarray}
S_2^2(\beta') &=& C_2^2(\alpha)S_2^2(\gamma) + S_2^2(\alpha)C_2^2(\gamma)
+ 2C_2(\alpha)S_2(\alpha)C_2(\gamma)S_2(\gamma)C_1(b)\nonumber\\
&+& k_1k_2S_2^2(\alpha)S_2^2(\gamma)S_1^2(b),\nonumber\\
T_2^2(\beta') &=& \frac{T_2^2(\alpha) + T_2^2(\gamma) + 2T_2(\alpha)T_2(\gamma)C_1(b)
+ k_1k_2T_2^2(\alpha)T_2^2(\gamma)S_1^2(b)}{(1 - k_2T_2(\alpha)T_2(\gamma)C_1(b))^2},\label{cos2b}
\end{eqnarray}
and
\begin{eqnarray}
S_2^2(\gamma) &=& C_2^2(\alpha)S_2^2(\beta') + S_2^2(\alpha)C_2^2(\beta')
- 2C_2(\alpha)S_2(\alpha)C_2(\beta')S_2(\beta')C_1(c)\nonumber\\
&+& k_1k_2S_2^2(\alpha)S_2^2(\beta')S_1^2(c),\nonumber\\
T_2^2(\gamma) &=& \frac{T_2^2(\alpha) + T_2^2(\beta') - 2T_2(\alpha)T_2(\beta')C_1(c)
+ k_1k_2T_2^2(\alpha)T_2^2(\beta')S_1^2(c)}{(1 + k_2T_2(\alpha)T_2(\beta')C_1(c))^2}.\label{cos2c}
\end{eqnarray}

What does it mean for triangle? From the Sine law (\ref{sinlaw}), having function $S(x)$ monotonically
increasing result that the longest side of any triangle is opposite to the largest angle and the
shortest side is opposed to the smallest angle. From the Cosine I law in its form that uses $C_1(x)$
function, having
$$C_1(b)C_1(c) + k_1S_1(b)S_1(c) = C_1(b-c),$$
$$C_1(a)C_1(b) - k_1S_1(a)S_1(b) = C_1(a+b),$$
and $C_i(x) le 1$ and is decreasing when $k_i = 1$, $C_i(x) = 1$ and is constant when $k_i = 0$,
$C_i(x) ge 1$ is increasing when $k_i = -1$, we can see:
$$C_1(a) = C_1(b)C_1(c) + k_1S_1(b)S_1(c)C_2(\alpha)$$
is equivalent to
$$a\begin{cases}
> b-c,\, k_2 = 1\\
= b-c,\, k_2 = 0\\
< b-c,\, k_2 = -1
\end{cases}$$
Similarly,
$$C_1(b) = C_1(a)C_1(c) - k_1S_1(a)S_1(c)C_2(\beta')$$
is equivalent to
$$b\begin{cases}
< a+c,\, k_2 = 1\\
= a+c,\, k_2 = 0\\
> a+c,\, k_2 = -1
\end{cases}$$

From the Cosine II law we can see:
$$C_2(\alpha) = C_2(\beta')C_2(\gamma) + k_2S_2(\beta')S_2(\gamma)C_1(a)$$
is equivalent to
$$\alpha\begin{cases}
> \beta'-\gamma,\, k_1 = 1\\
= \beta'-\gamma,\, k_1 = 0\\
< \beta'-\gamma,\, k_1 = -1
\end{cases}$$
Similarly,
$$C_2(\beta') = C_2(\alpha)C_2(\gamma) - k_2S_2(\alpha)S_2(\gamma)C_1(b)$$
is equivalent to
$$\beta'\begin{cases}
< \alpha+\gamma,\, k_1 = 1\\
= \alpha+\gamma,\, k_1 = 0\\
> \alpha+\gamma,\, k_1 = -1
\end{cases}$$

\section{Right (Quasi)--Triangle Equations}
We can define orthogonality in $\mathbb{B}^n$ using the orthogonality in $\mathbb{RP}^n$. Namely, two
vectors $v_1$ and $v_2$ of space $\mathbb{RP}^n$ are orthogonal, if $v_1\odot v_2 = 0$.

For $\mathbb{B}^2$ plane and a line, the orthogonal bundle is line only if $k_2 = 1$. In this case
when line rotates count--clockwise, its orthogonal line rotates count--clockwise and vice--versa (Figure
\ref{bundles} a). When $k_2 = 0$ there is the only orthogonal bundle, which doesn't rotate (Figure
\ref{bundles} b). When $k_2 = -1$ the orthogonal bundle rotates clockwise when the line rotates
count--clockwise toward to the same limit bundle and vice--versa (Figure \ref{bundles} c).

Generally we can't speak about right triangle as one of its catheti is line and another isn't (when
$k_2 \ne 1$). However, as we will see, this figure is important. We will name it right (quasi)--triangle,
which means right triangle, when $k_2 = 1$ and right quasi--triangle when $k_2 \ne 1$.

\begin{figure}[h]
\begin{center}
\includegraphics{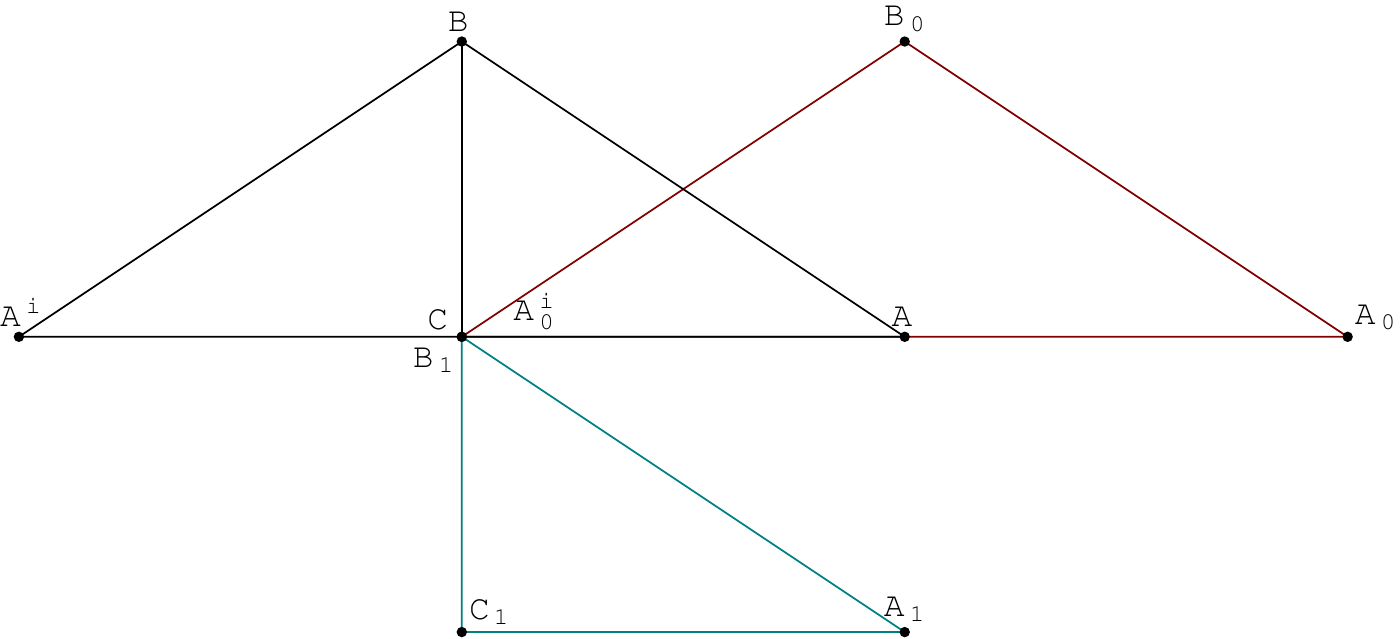}
\caption{Right (quasi)--triangle equations deduction.}
\label{right.triangle}
\end{center}
\end{figure}

We will construct a (quasi)--triangle as half of isosceles one (Figure \ref{right.triangle}). Consider a
triangle $A_0B_0A'_0$ with $A_0B_0 = A'_0B_0 = c$, $A_0A'_0 = 2b$, $\angle A'_0A_0B_0 = \angle A_0A'_0B_0
= \alpha$ and external angle $\angle A_0B_0A'_0 = 2\beta'$.

Let $A'_0 = O = \left[1: 0: 0\right]$ be origin, $A_0 = \mathfrak{R}_1(2b) A'_0 = \left[C_1(2b):
S_1(2b): 0\right]$, $B_0 = \mathfrak{R}_2(\alpha) \mathfrak{R}_1(c)$ $A_0 = \left[C_1(c): S_1(c)
C_2(\alpha): S_1(c) S_2(\alpha)\right]$.

Let $ABA' = \mathfrak{R}_1(-b)(A_0B_0A'_0)$ (Figure \ref{right.triangle}, black). Now, $A' = \mathfrak{R}_1
(-b)A'_0 = \left[C_1(b): -S_1(b):\right.$ $\left.0\right]$, $A = \mathfrak{R}_1(-b)A_0 = \left[C_1(b):
S_1(b):0\right]$ and $B = \mathfrak{R}_1(-b)B_0 =$
$$\begin{pmatrix}
C_1(b)&k_1S_1(b)&0\\
-S_1(b)&C_1(b)&0\\
0&0&1
\end{pmatrix}
\begin{pmatrix}
C_1(c)\\
S_1(c)C_2(\alpha)\\
S_1(c)S_2(\alpha)
\end{pmatrix}$$
$= \left[C_1(b)C_1(c) + k_1S_1(b)S_1(c)C_2(\alpha):-S_1(b)C_1(c) + C_1(b)S_1(c)C_2(\alpha):
S_1(c)\right.$ $\left.S_2(\alpha)\right]$.

Finally, let $C \in AA', AC = A'C = b$. Then $C = \left[1:0:0\right]$ is origin. From figure equality
$A'BC = ABC$ result $BC \perp A'A$. Therefore we can consider $ABC$ right (quasi)--triangle.

Having figures $A'BC = ABC$ and $C$ is origin, result $B$ have form $B = (x, 0, y)$, where
\begin{eqnarray*}
C_1(b)C_1(c) + k_1S_1(b)S_1(c)C_2(\alpha) &=& x,\\
-S_1(b)C_1(c) + C_1(b)S_1(c)C_2(\alpha) &=& 0,\\
S_1(c)S_2(\alpha) &=& y.
\end{eqnarray*}

From the second equality, have
\begin{equation}\label{nearangle1}
T_1(b) = T_1(c)C_2(\alpha).
\end{equation}
Using the value of $C_2(\alpha)$ from this equality and putting it in the first one, have
\begin{eqnarray*}
x &=& C_1(b) C_1(c) + k_1S_1(b)S_1(c) \frac{T_1(b)}{T_1(c)}\\
&=& \frac{C_1(c)}{C_1(b)}(C_1^2(b) + k_1S_1^2(b)) = \frac{C_1(c)}{C_1(b)}.
\end{eqnarray*}

Let now calculate the value of
\begin{eqnarray*}
x^2 + k_1k_2y^2 &=& \frac{C_1^2(c)}{C_1^2(b)} + k_1k_2S_1^2(c)S_2^2(\alpha)\\
&=& \frac{C_1^2(c)}{C_1^2(b)} + k_1S_1^2(c)(1 - C_2^2(\alpha))\\
&=& \frac{C_1^2(c)}{C_1^2(b)} + k_1S_1^2(c) - k_1S_1^2(c)\frac{T_1^2(b)}{T_1^2(c)}\\
&=& \frac{C_1^2(c)}{C_1^2(b)} + k_1S_1^2(c) - k_1S_1^2(b)\frac{C_1^2(c)}{C_1^2(b)}\\
&=& \frac{C_1^2(c)}{C_1^2(b)}(1 - k_1S_1^2(b)) + k_1S_1^2(c)\\
&=& \frac{C_1^2(c)}{C_1^2(b)}C_1^2(b) + k_1S_1^2(c)\\
&=& C_1^2(c) + k_1S_1^2(c) = 1.
\end{eqnarray*}
It means that exists $a \in \mathbb{R}, C_{12}(a) = \frac{C_1(c)}{C_1(b)}, S_{12}(a) = S_1(c)S_2
(\alpha)$ that has characteristic $k = k_1k_2$. It is a `distance' parameter $BC$. We have two more
equations:
\begin{eqnarray}
C_1(c) &=& C_{12}(a)C_1(b)\label{pitagora}\\
S_{12}(a) &=& S_1(c)S_2(\alpha)\label{awayangle1}
\end{eqnarray}
From (\ref{nearangle1}) and (\ref{awayangle1}) have:
$$\frac{S_{12}(a)}{T_1(b)} = \frac{S_1(c)S_2(\alpha)}{T_1(c)C_2(\alpha)} = C_1(c)T_2(\alpha),$$
using the value of $C_1(c)$ from (\ref{pitagora}) have
\begin{eqnarray}
\frac{S_{12}(a)}{T_1(b)} &=& C_{12}(a)C_1(b)T_2(\alpha),\nonumber\\
T_{12}(a) &=& S_1(b)T_2(\alpha).\label{catheti}
\end{eqnarray}

The last 6 equations will include $\beta'$. In order to be able to deduce them we will introduce
translation $\mathfrak{T}(-a)$, so as $\mathfrak{T}(-a)B = C$. Having characteristic $a$ is $k_1k_2
= K_2$,
$$\mathfrak{T}(-a) =
\begin{pmatrix}
C_{12}(a)&0&K_2S_{12}(a)\\
0&1&0\\
-S_{12}(a)&0&C_{12}(a)
\end{pmatrix}.$$
We can check this map preserves vector product.

Applying $\mathfrak{T}(-a)$, obtain $B_1 = \mathfrak{T}(-a)B = \left[1:0:0\right]$, $C_1 =
\mathfrak{T}(-a)C = \left[C_{12}(a):0:-S_{12}(a)\right]$ and $A_1 = \mathfrak{T}(-a)A =$
$$\begin{pmatrix}
C_{12}(a)&0&K_2S_{12}(a)\\
0&1&0\\
-S_{12}(a)&0&C_{12}(a)
\end{pmatrix}
\begin{pmatrix}
C_1(b)\\
S_1(b)\\
0
\end{pmatrix}$$
$= \left[C_{12}(a)C_1(b):S_1(b):-S_{12}(a)C_1(b)\right] = \left[C_1(c): S_1(c)C_2(\beta'):
-S_1(c)S_2(\beta')\right]$ (Figure \ref{right.triangle}, cyan). From here we have
\begin{equation}\label{awayangle2}
S_1(b) = S_1(c)C_2(\beta').
\end{equation}
Moreover, having (\ref{pitagora}), obtain
\begin{eqnarray}
S_{12}(a)C_1(b) &=& S_1(c)S_2(\beta'),\nonumber\\
S_{12}(a)\frac{C_1(c)}{C_{12}(a)} &=& S_1(c)S_2(\beta'),\nonumber\\
T_{12}(a) &=& T_1(c)S_2(\beta').\label{nearangle2}
\end{eqnarray}
Combining the last 2 equalities (\ref{awayangle2}), (\ref{nearangle2}) with (\ref{pitagora}), we have
\begin{eqnarray}
\frac{T_{12}(a)}{S_1(b)} &=& \frac{T_1(c)S_2(\beta')}{S_1(c)C_2(\beta')}\nonumber\\
= \frac{T_2(\beta')}{C_1(c)} &=& \frac{T_2(\beta')}{C_{12}(a)C_1(b)},\nonumber\\
S_{12}(a) &=& T_1(b)T_2(\beta').\label{catheti2}
\end{eqnarray}

Now, having (\ref{nearangle1}) and (\ref{awayangle2}):
$$T_1(c)C_2(\alpha) = T_1(b) = \frac{S_1(b)}{C_1(b)} = \frac{S_1(c)C_2(\beta')}{C_1(b)},$$
calculate with (\ref{pitagora}):
\begin{eqnarray}
C_1(b) &=& \frac{S_1(c)C_2(\beta')}{T_1(c)C_2(\alpha)}\nonumber\\
&=& C_1(c)\frac{C_2(\beta')}{C_2(\alpha)}\nonumber\\
&=& C_{12}(a)C_1(b)\frac{C_2(\beta')}{C_2(\alpha)},\nonumber\\
C_{12}(a)\frac{C_2(\beta')}{C_2(\alpha)} &=& 1,\nonumber\\
C_2(\alpha) &=& C_{12}(a)C_2(\beta').\label{angles3}
\end{eqnarray}

Now from (\ref{catheti}), (\ref{awayangle2}) and (\ref{nearangle2}),
\begin{eqnarray}
T_{12}(a) = S_1(b)T_2(\alpha) &=& T_1(c)S_2(\beta'),\nonumber\\
S_1(c)C_2(\beta')T_2(\alpha) &=& T_1(c)S_2(\beta'),\nonumber\\
T_2(\beta') &=& C_1(c)T_2(\alpha).\label{angles1}
\end{eqnarray}

Finally, by multiplying the last equations (\ref{angles3}) and (\ref{angles1}), have using
(\ref{pitagora}):
\begin{eqnarray}
C_1(c)T_2(\alpha)C_2(\alpha) &=& T_2(\beta')C_{12}(a)C_2(\beta'),\nonumber\\
C_1(c)S_2(\alpha) &=& C_{12}(a)C_2(\beta'),\nonumber\\
C_{12}(a)C_1(b)S_2(\alpha) &=& C_{12}(a)S_2(\beta'),\nonumber\\
S_2(\beta') &=& C_1(b)S_2(\alpha).\label{angles2}
\end{eqnarray}

It is necessary to modify equations (\ref{pitagora}) and (\ref{angles3}) in order to not contain the
$C(x)$ function.
\begin{eqnarray*}
k_1S_1^2(c) &=& 1 - C_1^2(c) = (C^2_{12}(a) + k_1k_2S^2_{12}(a))(C_1^2(b) + k_1S_1^2(b))
- C^2_{12}(a)C_1^2(b)\\
&=& k_1k_2S^2_{12}(a)C_1^2(b) + k_1C^2_{12}(a)S_1^2(b) + k_1^2k_2S^2_{12}(a)S_1^2(b),\\
S_1^2(c) &=& k_2S^2_{12}(a)C_1^2(b) + C^2_{12}(a)S_1^2(b) + k_1k_2S^2_{12}(a)S_1^2(b)
\end{eqnarray*}
By dividing the last equality by its $C(x)$ form, obtain:
\begin{equation}\label{pitagora1}
T_1^2(c) = k_2T^2_{12}(a) + T_1^2(b) + k_1k_2T^2_{12}(a)T_1^2(b).
\end{equation}
Similarly,
\begin{eqnarray*}
k_2S_2^2(\alpha) &=& 1 - C_2^2(\alpha)
= (C^2_{12}(a) + k_1k_2S^2_{12}(a))(C_2^2(\beta') + k_2S_2^2(\beta')) - C^2_{12}(a)C_2^2(\beta')\\
&=& k_2C^2_{12}(a)S_2^2(\beta') + k_1k_2S^2_{12}(a)C_2^2(\beta') + k_1k_2^2S^2_{12}(a)S_2^2(\beta'),\\
S_2^2(\alpha) &=& C^2_{12}(a)S_2^2(\beta') + k_1S^2_{12}(a)C_2^2(\beta')
+ k_1k_2S^2_{12}(a)S_2^2(\beta')
\end{eqnarray*}
By dividing the last equality by its $C(x)$ form, obtain:
\begin{equation}\label{angles4}
T_2^2(\alpha) = k_1T^2_{12}(a) + T_2^2(\beta') + k_1k_2T^2_{12}(a)T_2^2(\beta')
\end{equation}

Note that for $k_2 = 1$ equations (\ref{nearangle1}) --- (\ref{angles4}) can be used if external
angle $\beta'$ change to internal $\beta$ with the following changes:
\begin{eqnarray*}
\beta &=& \frac{\pi}2 - \beta'\\
\cos\beta &=& \sin\beta'\\
\sin\beta &=& \cos\beta'\\
\tan\beta &=& \cot\beta'\\
\cot\beta &=& \tan\beta'
\end{eqnarray*}

\section{More rotations}
As we can see, transformation $\mathfrak{T}(-a)$ preserves vector product. In order to be a
motion it needs to be presented as finite product of main rotations. If $a$, $b$, $c$, $\alpha$
and $\beta'$ are real numbers for which have place equalities (\ref{nearangle1}) --- (\ref{angles2})
then it can be checked that $\mathfrak{T}(a) = \mathfrak{R}_2(\beta')\mathfrak{R}_1(c)\mathfrak{R}_2
(-\alpha)\mathfrak{R}_1(-b)$. We will introduce new transformations as following:
$$\mathfrak{R}_{ij}(\phi) =
\begin{pmatrix}
1&\ldots&0&\ldots&0&\ldots&0\\
\vdots&\ddots&\vdots&\ddots&\vdots&\ddots&\vdots\\
0&\ldots&C_{i+1,...j}(\phi)&\ldots&-\frac{K_j}{K_i}S_{i+1,...j}(\phi)&\ldots&0\\
\vdots&\ddots&\vdots&\ddots&\vdots&\ddots&\vdots\\
0&\ldots&S_{i+1,...j}(\phi)&\ldots&C_{i+1,...j}(\phi)&\ldots&0\\
\vdots&\ddots&\vdots&\ddots&\vdots&\ddots&\vdots\\
0&\ldots&0&\ldots&0&\ldots&1
\end{pmatrix}.$$

It's easy to see that all $\mathfrak{R}_{ij}(\phi)$ are motions. All they can be presented
as finit product of main rotations. We will name them \textit{rotations} of the space $\mathbb{B}^n$.
In special case, $\mathfrak{R}_{0i}(\phi)$ we will name them \textit{translations} $\mathfrak{T}_{i}
(\phi)$ of the space.

\section{Generalized Orthogonal Matrix}
For $\mathbb{RP}^n$ with given specification $k_p,\, p=\overline{1,n}$ we will name the vector
$x \in \mathbb{RP}^n$ \textit{upper $i$-normalized}, $i \in \overline{0,n}$ if $\frac1{K_i}x\odot
x = 1$. For $0 \le i < j \le n$ we will name two vectors $x$ and $y$ \textit{upper $ij$-orthogonal} if
$\frac1{K_{min(i,j)}}x\odot y = 0$. We will name the matrix $M_{(n+1)\times(n+1)}$ composed of
columns $c_i$ \textit{upper orthogonal} if all columns $c_i$ are upper $i$-normalized and any two
columns $c_i$ and $c_j$ are upper $ij$-orthogonal.

It's easy to see that all main rotation matrixes are upper orthogonal. Moreover product of two upper
orthogonal matrix is upper orthogonal. Really, let $X,Y$ are two upper orthogonal matrices. It means
that $X$ is composed of $(x_0,...,x_n)$ columns and $Y$ --- from $(y_0,...,y_n)$ columns, and
$\frac1{K_{min(i,j)}}x_i\odot x_j = \frac1{K_{min(i,j)}}y_i\odot y_j = \delta_{ij}$ for all $i,j =
\overline{0,n}$, where $\delta_{ij} = 1, i=j$ and $\delta_{ij} = 0, i\ne j$. Let $Z=XY$ with
elements $z_{ij} = \sum_{p=0}^nx_{ip}y_{pj}$. Let $z_i$ and $z_j$ be 2 columns of $Z$. Let calculate
\begin{eqnarray*}
\frac1{K_{min(i,j)}}z_i\odot z_j &=& \frac1{K_{min(i,j)}}\sum_{p=0}^nK_pz_{pi}z_{pj}\\
&=& \frac1{K_{min(i,j)}}\sum_{p=0}^nK_p\left(\sum_{m_1=0}^nx_{pm_1}y_{m_1i}\right)
\left(\sum_{m_2=0}^nx_{pm_2}y_{m_2j}\right)\\
&=& \frac1{K_{min(i,j)}}\sum_{p=0}^nK_p\sum_{m_1=0}^n\sum_{m_2=0}^nx_{pm_1}x_{pm_2}y_{m_1i}y_{m_2j}\\
&=& \frac1{K_{min(i,j)}}\sum_{m_1=0}^n\sum_{m_2=0}^ny_{m_1i}y_{m_2j}\sum_{p=0}^nK_px_{pm_1}x_{pm_2}\\
&=& \frac1{K_{min(i,j)}}\sum_{m_1=0}^n\sum_{m_2=0}^ny_{m_1i}y_{m_2j}K_{min(m_1,m_2)}\delta_{m_1m_2}\\
&=& \frac1{K_{min(i,j)}}\sum_{m=0}^ny_{mi}y_{mj}K_m\\
&=& \delta_{ij}
\end{eqnarray*}

We will name the vector $x \in \mathbb{RP}^n$ \textit{lower $i$-normalized}, $i \in \overline{0,n}$
if $K_i\sum_{j=0}^n\frac{x_j^2}{K_j} = 1$. For $0 \le i < j \le n$ we will name two vectors $x$ and
$y$ \textit{lower $ij$-orthogonal} if $K_{max(i,j)}\sum_{p=0}^n\frac{x_py_p}{K_p} = 0$. We will name
the matrix $M_{(n+1)\times(n+1)}$ composed of rows $r_i$ \textit{lower orthogonal} if all rows $r_i$
are lower $i$-normalized and any two rows $r_i$ and $r_j$ are lower $ij$-orthogonal.

It's easy to see that all main rotation matrixes are also lower orthogonal. Moreover product of two
lower orthogonal matrices is lower orthogonal. Really, let $X,Y$ are two lower orthogonal matrices.
It means that $X$ is composed of $(x_0,...,x_n)$ rows and $Y$ --- from $(y_0,...,y_n)$ rows, where
$K_{max(i,j)}\sum_{p=0}^n\frac{x_{ip}x_{jp}}{K_p} = K_{max(i,j)}\sum_{p=0}^n\frac{y_{ip}y_{jp}}{K_p}
= \delta_{ij}$ for all $i,j = \overline{0,n}$. Let $Z=XY$. Let calculate
\begin{eqnarray*}
K_{max(i,j)}\sum_{p=0}^n\frac{z_{ip}z_{jp}}{K_p}
&=& K_{max(i,j)}\sum_{p=0}^n\frac1{{K_p}}\left(\sum_{m_1=0}^nx_{im_1}y_{m_1p}\right)
\left(\sum_{m_2=0}^nx_{jm_2}y_{m_2p}\right)\\
&=& K_{max(i,j)}\sum_{p=0}^n\frac1{K_p}\sum_{m_1=0}^n\sum_{m_2=0}^nx_{im_1}x_{jm_2}y_{m_1p}y_{m_2p}\\
&=& K_{max(i,j)}\sum_{m_1=0}^n\sum_{m_2=0}^nx_{im_1}x_{jm_2}\sum_{p=0}^n\frac1{K_p}y_{m_1p}y_{m_2p}\\
&=& K_{max(i,j)}\sum_{m_1=0}^n\sum_{m_2=0}^nx_{im_1}x_{jm_2}\frac{\delta_{m_1m_2}}{K_{max(m_1,m_2)}}\\
&=& K_{max(i,j)}\sum_{m=0}^n\frac{x_{im}x_{jm}}{K_m}\\
&=& \delta_{ij}
\end{eqnarray*}

For some upper orthogonal matrix $X$ has place the equality
\begin{eqnarray*}
\frac1{K_j}\sum_{i=0}^nK_ix^2_{ij} &=& 1,\\
\sum_{i=0}^nK_ix^2_{ij} &=& K_j,\\
\sum_{i=0}^nx^2_{ij}\prod_{p=1}^ik_p &=& \prod_{p=1}^jk_p.
\end{eqnarray*}
Let divide it to $K_q$, $q \le j$:
$$\sum_{i=0}^{q-1}\frac{K_i}{K_q}x^2_{ij} + \sum_{i=q}^n\frac{K_i}{K_q}x^2_{ij} = \frac{K_j}{K_q}.$$
As $K_q$ divides $K_j$ and $K_i, i \ge q$, but $K_q$ doesn't divide $K_i, i < q$, result that for
$0 \le i < j \le n$, $x_{ij}$ divide $K_q/K_i$ for all $i < q \le j$, or $x_{ij}$ divide $K_j/K_i$.

Having for some upper orthogonal matrix $X$, elements $x_{ij}$ divide $K_j/K_i$, construct the matrix
$Y$ of the same size with elements $y_{ij} = \sqrt{\frac{K_i}{K_j}}x_{ij}$. The matrix $Y$ is
orthogonal one\footnote{consider $\sqrt{0} = 0$.} (if may be complex, in this case it isn't unitar, but
orthogonal). Really, for some $i = \overline{0, n}$,
$$\sum_{i=0}^ny^2_{ij} = \sum_{i=0}^n\frac{K_i}{K_j}x^2_{ij} = \frac1{K_j}\sum_{i=0}^nK_ix^2_{ij} = 1$$
and for some $j_1 \ne j_2 = \overline{0,n}$,
\begin{eqnarray*}
\sum_{i=0}^ny_{ij_1}y_{ij_2} &=& \sum_{i=0}^n\frac{K_i}{\sqrt{K_{j_1}K_{j_2}}}x_{ij_1}x_{ij_2}\\
&=& \frac{K_{min(j_1,j_2)}}{\sqrt{K_{j_1}K_{j_2}}}\frac1{K_{min(j_1,j_2)}}
\sum_{i=0}^nK_ix_{ij_1}x_{ij_2}\\
&=& \frac{K_{min(j_1,j_2)}}{\sqrt{K_{j_1}K_{j_2}}}0 = 0,
\end{eqnarray*}
because $x_{ij_1}x_{ij_2}$ divides $\frac{\sqrt{K_{j_1}K_{j_2}}}{K_i}$.

For some orthogonal matrix always has place also the following equalities for $i_1 \ne i_2 =
\overline{0,n}$:
$$\sum_{j=0}^ny^2_{ij} = \sum_{j=0}^n\frac{K_i}{K_j}x^2_{ij} = K_i\sum_{j=0}^n\frac1{K_j}x^2_{ij} = 1,$$
\begin{eqnarray*}
\sum_{j=0}^ny_{i_1j}y_{i_2j}
&=& \sum_{j=0}^n\frac{\sqrt{K_{i_1}K_{i_2}}}{K_j}x_{i_1j}x_{i_2j}\\
&=& \frac{\sqrt{K_{i_1}K_{i_2}}}{K_{max(i_1,i_2)}}K_{max(i_1,i_2)}
\sum_{j=0}^n\frac1{K_j}x_{i_1j}x_{i_2j}\\
&=& \frac{\sqrt{K_{i_1}K_{i_2}}}{K_{max(i_1,i_2)}}0 = 0
\end{eqnarray*}
So, $x_{i_1j}x_{i_2j}$ divides $\frac{K_j}{\sqrt{K_{i_1}K_{i_2}}}$. It means that $X$ is also lower
orthogonal matrix.

Inverse orthogonal matrix $Y^{-1}$ is easy constructed as $y'_{ji} = y_{ij}$. Then
\begin{eqnarray*}
\sqrt{\frac{K_j}{K_i}}x'_{ji} &=& \sqrt{\frac{K_i}{K_j}}x_{ij},\\
x'_{ji} &=& \frac{K_i}{K_j}x_{ij}.
\end{eqnarray*}
The last equality isn't applicable if some characteristic $k_i = 0$. Although is true, it isn't
determinable having the form of $0/0$. If some characteristic $k_m = 0, m < n$, the matrix has the form:
$$M = \begin{pmatrix}
A&O\\
B&C
\end{pmatrix}$$
Really, for the first $m$ columns the upper orthogonality condition is equivalent to:
\begin{eqnarray*}
\sum_{i=0}^n\frac{K_i}{K_j}x_{ij}^2 &=& \sum_{i=0}^{m-1}\frac{K_i}{K_j}x_{ij}^2, \forall
j=\overline{0,m-1},\\
\sum_{i=0}^n\frac{K_i}{K_{j_1}}x_{ij_1}x_{ij_2} &=& \sum_{i=0}^{m-1}\frac{K_i}{K_{j_1}}x_{ij_1}x_{ij_2},
\forall j_1 =\overline{0,m-1}, j_2=\overline{0,n}, j_1 < j_2,
\end{eqnarray*}
Having $K_j \ne 0$ and $K_i = 0$, all terms, starting with $i = m$ equals to zero. So, matrix $A$ is
upper orthogonal of size $m\times m$ and matrix $B$ is free of size $(n-m+1)\times m$. For the last
$n-m+1$ columns upper orthogonality has form:
\begin{eqnarray*}
\sum_{i=0}^n\frac{K_i}{K_j}x_{ij}^2 &=& \sum_{i=m}^{n}\frac{\prod_{p=m}^ik_p}{\prod_{p=m}^jk_p}x_{ij}^2,
\forall j=\overline{m,n},\\
\sum_{i=0}^n\frac{K_i}{K_{j_1}}x_{ij_1}x_{ij_2} &=&
\sum_{i=m}^{n}\frac{\prod_{p=m}^ik_p}{\prod_{p=m}^{j_1}k_p}x_{ij_1}x_{ij_2}, \forall j_1 < j_2 =
\overline{m,n},
\end{eqnarray*}
because $K_j = 0$. It means the matrix $C$ is upper orthogonal of size $(n-m+1)\times(n-m+1)$ and the
matrix $O$ is obligatory zero one of size $m\times(n-m+1)$ (otherwise elements of $M$ aren't finite).

It's easy to verify that inverse matrix has form:
$$M^{-1} = \begin{pmatrix}
A^{-1}&O\\
-C^{-1}BA^{-1}&C^{-1}\end{pmatrix}$$
This way of calculating the inverse matrix can easy be generalized to either number of null
characteristics.

We will name upper orthogonal matrixes (which also are lower orthogonal) \textit{generalized orthogonal}.
We will use the term orthogonal matrix meaning generalized orthogonal matrix if isn't stated otherwise.
As we can see, the orthogonal matrix set is closed in respect of multiplying, it contains the unit
element and for any element it contains its inverse. So the orthogonal matrix set form isomentry group
of space. All motion matrices are generalized orthogonal.

\section[Orthogonal Matrix Decomposition]{Orthogonal Matrix as Product of Rotations}
Let $X$ will be orthogonal matrix. We will search the rotation matrices, the product of which gives $X$.
Note that The matrix $X\mathfrak{R}_{ij}(\phi)$ have all columns $x_p$ of $X$ except $i$ and $j$ ones.
These columns are $x'_i = x_iC_{i+1...j}(\phi) + x_jS_{i+1...j}(\phi)$ and $x'_j = -\frac{K_j}{K_i}x_i
S_{i+1...j}(\phi) + x_jC_{i+1...j}$.

For the last row let separate elements $x_i, i = \overline{0,n}$ in three categories: having
characteristics $K_n/K_i$ equals to 1, 0 and $-1$. Note that for $i$ row the $i$ element is always of the
category 1, because its characteristic is $K_i/K_i = 1$. We will multiply $X$ on the right by
$\mathfrak{R}_{in}(\phi), i = \overline{0,n}$ in order to have in the $n$-th row a single element of
category 1 and single element of category $-1$, different from 0. All these rotations are elliptic ones.
For elements of the same characteristic $x_{ni}$ and $x_{nj}$ we can use $\cos\phi = \frac{x_{ni}}
{\sqrt{x^2_{ni}+x^2_{nj}}}$ and $\sin\phi = -\frac{x_{nj}}{\sqrt{x^2_{ni}+x^2_{nj}}}$. Moreover,
always $x_{nn} \ne 0$.

Now, we have one element of category 1 and one of category $-1$, different from zero (the $n$-th one
and, for example, the $p$-th one) and element of category 1 has absolute value greater then the element
of category $-1$ because for this lower $n$-normalized row have place equality $x^2_{nn} - x^2_{np} = 1$.
It means that exists $\phi \in \mathbb{R}$ so that $\cosh\phi = \frac{x_{nn}}{\sqrt{x^2_{nn} -
x^2_{np}}}$ and $\sinh\phi = -\frac{x_{np}}{\sqrt{x^2_{nn} - x^2_{np}}}$ and hyperbolic rotation
that transforms the element of category $-1$, $x_{np}$ in $x'_{np} = 0$ and the element of category 1,
$x_{nn}$ in $x'_{nn} \ne 0$.

For category 0 there exist parabolic rotations, that preserves the element of category 1 ($x_{nn}$) and
elements of category 0 transform in 0. For this case, if one this element is on $q$-th column, $\phi =
-x_{nq}$. The last non--zero element $x_{nn}$ equals to 1 or $-1$, because the last row is lower
$n$-normalized.

We can consider the first $n$ columns as having $n$ elements (the last one equals to zero). They form
orthogonal matrix of size $n$. The last, $(n+1)$-th column (without the last element) is upper
$in$-orthogonal to first $n$ columns, $i = \overline{0, n-1}$. As these columns have $n$ elements each,
$(n+1)$-th column is obligatory null (excluding the last its element).

In this stage we can consider the resulting matrix as having size $n$ instead of $n+1$ and repeat the
process for it. Finally obtain the matrix $E$ which has elements on main diagonal 1 or $-1$ and all rest
elements 0. It is the reflection matrix on a point or line or plane or hyperplane. Obtain the equality:
$X\prod_{j=1}^q\mathfrak{M}_j = E$. It's easy to see that $X = E\prod_{j=q}^1\mathfrak{M}^{-1}_j$
($q = (n+1)n/2$). Strictly speaking, the matrix $E$ can't be presented as product of rotations. In order
to identify motions of $\mathbb{B}^n$ with orthogonal matrices, we should name $E$ (which preserve
vector product) motion. However, these motions are improper (there is no continuous
parameterization of motion $\mathfrak{M}(\alpha)$ on a segment $[0,1]$ such that $\mathfrak{M}(0) =
I, \mathfrak{M}(1) = E$ and all $\mathfrak{M}(\alpha)$ are motions on all $\alpha\in[0,1]$).
Having in expression for $X$ determinant of matrices $\mathfrak{M}_i$ equals to 1 and determinant of
$E$ is $\pm 1$, determinant of $X$ equals $\pm 1$.

\section{Coordinate and State Matrix}
Consider in some space $\mathbb{B}^n$ $n + 1$ vectors $v_0,... v_n$. Let coordinates of $v_i$ are
$\left[v_{0i}:...:v_{ni}\right], i = \overline{0,n}$. Let vectors $v_i$ are ordered and form basis
of $\mathbb{B}^n$ (not obligatory orthonormal). Compose the matrix $V$ with elements $v_{ij}, i,j =
\overline{0,n}$. Will name it \textit{coordinate matrix} for vectors $v_i$. Construct also the matrix
$M$ of size $(n + 1)\times(n + 1)$ with elements $m_{i,j} = \frac{v_i\odot v_j}{K_i}$. We will name
the matrix $M$ \textit{state matrix} of $v_i$. Having $v_i$ is space basis, elements $m_{ij}$ are all
finite. State matrix shows how orthonormal is some vector family. It tends to unite one when vectors
are more normalized and orthogonal to each other\footnote{If space specification contains null
characteristics, some elements of state matrix can have any value, even for orthonormal vector family.}.
We will demonstrate that volume of parallelepiled constructed on vectors $v_i$ equals to $|\det V|$ and
\begin{equation}\label{det}
\det M = \left(\det V\right)^2.
\end{equation}

First, let $v_i, i = \overline{0,n}$ are orthonormal. Then the parallelepiped volume is 1, the matrix
$V$ is orthogonal one and $\det V = \pm 1$. So, $|\det V|$ equals to parallelepiped volume. All elements
on main diagonal $m_{ii} = 1$, because all vectors $v_i$ are upper $i$-normalized. All elements above main
diagonal are $m_{ij} = 0, i < j$, because all vectors $v_i$ and $v_j$ are upper $ij$-orthogonal (elements
under the main diagonal may differ from 0). It means that the matrix $M$ is lower triangular with all
elements on main diagonal equls to 1 and $\det M = 1 = \left(\det V\right)^2$.

Further, note that $x\odot(y + z) = x\odot y + x\odot z, (x + y)\odot z = x\odot z + y\odot z,
(\alpha x)\odot y = \alpha(x\odot y) = x\odot(\alpha y), \forall x, y, z \in\mathbb{B}^n, \alpha
\in\mathbb{R}$. Matrix determinant equals to zero if it contain proportional columns or rows. When some
row or column of a matrix is multiplied by $\alpha$, the resulting matrix determinant is $\alpha$ times
original matrix determinant. When some row or column is sum of two rows / columns, then the matrix
determinant equals to sum of determinants of matrices containing the first and the second row / column.

Second, let instead of some $v_i$ use $v'_i = \alpha v_i$. In this case parallelepiped volume grows in
$\alpha$ times and $|\det V'| = |\alpha||\det V|$. Moreover, $\det M' = \alpha^2\det M = \left(\alpha
\det V\right)^2 = \left(\det V'\right)^2$.

Third, let instead of some $v_i$ use $v'_i = v_i + \alpha v_j$. In this case parallelepiped volume
remains unchanged, as well as determinant of $V$, and $\det M' = \det M = \left(\det V\right)^2 = \left(
\det V'\right)^2$.

Finally, observe, that all matrices $V$ result from orthogonal matrices using operations form the second
and the third step. It means the equation (\ref{det}) is true for all matrices and volume of
parallelepiped constructed on vectors $v_i$ equals to $|\det V|$.

When somebody calculates the parallelepiped volume it's usefull to use the state matrix. Its elements
don't change on motions and it is always square, even when the number of vectors is less then the space
dimension (the matrix $V$ isn't square in this case).

\section[Plane definition and Specification]{Plane definition and Specification. Lineals and their
Specification}
Having $\mathbb{B}^m \subset \mathbb{RP}^m$, all $m$-dimensional planes $\mathbb{B}^n$ ($m < n$)
lie in $\mathbb{RP}^n$ with one global condition for vectors $X \in\mathbb{RP}^m\subset\mathbb
{PR}^n$: $X\odot X = 1$. Leaving this condition (it doesn't change on any motion), we can consider
$m$-dimensional planes of $\mathbb{B}^n$ as $m$-dimensional planes of $\mathbb{RP}^n$.

By definition, $m$-dimensional plane $L^m$ results from subspace $\mathbb{B}^m$ on some motion.
Subspace $\mathbb{B}^m$ has the first $m+1$ columns of unite matrix with dimension $n+1$ as its basis.
Multiplying the basis matrix of $\mathbb{B}^m$ by some orthogonal matrix result basis matrix of $L^m$
as first $m+1$ columns of orthogonal matrix. Being a subspace, specification of $\mathbb{B}^m$ contain
the first $m$ characteristics of specification $\mathbb{B}^n$.

What happens if we take any $m+1$ columns of some orthogonal matrix as basis? Let column indices $i_0,
i_1,... i_m$ and $i_p = \overline{0,n}$. It's easy to see that motions that preserve this figure only
change these columns (interior figure motions) or change no these columns (motions of $\mathbb{B}^n$
that preserve all its points). Thus, figure characteristics $K_p' = K_{i_p}, p = \overline{0,m}$, or
$k_p' = K_p'/K_{p-1}' = K_{i_p}/K_{i_{p-1}} = \prod_{j=i_{p-1}+1}^{i_p}k_j, p = \overline{1,m}$ is
its specification. These figures generally speaking are not planes. We will name them \textit{lineals}.
We will name planes also lineals.

It may happen, that some lineal has $K_0' \ne 1$ (space and planes have it equals to 1). In this case
lineal may not intersect the space sphere and may not have image. We will name lineals that have no
image improper. Although they have no image, their properties help studying the space geometry.

One more interesting case is when the space specification has characteristic $-1$ and some lineal is
constructed on limit vectors for this characteristic. Such lineals can't be constructed from matrices
get as finite product of motions. They can be constructed as limit of infinite products. These lineal
specifications can't be deduced from space specification.

For example, let space $\mathbb{B}^2$ has specification $\{1, -1\}$. Vectors $\left[0:1:1\right]$
and $\left[0:-1:1\right]$ can't result from coordinate vectors on finite product of motions. However,
there exist translations along these vectors (interior lineal translations):
$$\mathfrak{M} = \begin{pmatrix}
1&-2S_0(\frac{\phi}{2})&-2S_0(\frac{\phi}{2})\\
2S_0(\frac{\phi}{2})&1-2S_0^2(\frac{\phi}{2})&-2S_0^2(\frac{\phi}{2})\\
-2S_0(\frac{\phi}{2})&2S_0^2(\frac{\phi}{2})&1+2S_0^2(\frac{\phi}{2})
\end{pmatrix},$$
$$\mathfrak{W} = \begin{pmatrix}
1&2S_0(\frac{\psi}{2})&-2S_0(\frac{\psi}{2})\\
-2S_0(\frac{\psi}{2})&1-2S_0^2(\frac{\psi}{2})&2S_0^2(\frac{\psi}{2})\\
-2S_0(\frac{\psi}{2})&-2S_0^2(\frac{\psi}{2})&1+2S_0^2(\frac{\psi}{2})
\end{pmatrix}.$$
These motions matrices use functions $C_0(x)$ and $S_0(x)$ have characteristic 0, despite the fact
the space specification doesn't contain zero. These translations are border space motions between
elliptic and hyperbolic ones.

\section[Projection of Vector on Lineal]{Projection of Vector on Lineal and on its Orthogonal
Completion}
We will name some vector $v'$ projection of vector $v$ on lineal $L^m$, if $v' \in L^m, v - v' \perp
L^m$. Let lineal $L^m$ is constructed on vectors $l_0,... l_m$. Then $v' = \sum_{i=0}^m\frac{v\odot
l_i}{K_i}l_i$. Evident, $v' \in L_m$. Let's see:
\begin{eqnarray*}
(v - v')\odot l_j &=& \left(v - \sum_{i=0}^m\frac{v\odot l_i}{K_i}l_i\right)\odot l_j\\
&=& v\odot l_j - \sum_{i=0}^m\frac{v\odot l_i}{K_i}(l_i\odot l_j)\\
&=& v\odot l_j - \frac{v\odot l_j}{K_j}(l_j\odot l_j)\\
&=& v\odot l_j - \frac{v\odot l_j}{K_j}K_j = 0,
\end{eqnarray*}
for all $j = \overline{0,m}$, in other words, $v'' = v - v' \perp L^m$. When some $K_i = 0$,
expression $\frac{v\odot l_i}{K_i}$ has undefined value. It happens when some vector direction is
orthogonal to all others. In this case there is impossible to determine unique orthogonal vector.
However, any value of this expression, for example 0, is valid, as it corresponds to some orthogonal
vector.

\section[Basis Change in Lineal]{Basis Change in Lineal. Unique Form of Lineal}
Let $L^m \subset \mathbb{B}^n$ is some space lineal, defined by matrix of size $(n+1)\times(m+1)$. The
matrix columns $l_i, i=\overline{0,m}$ form basis of lineal. Consider vector $a = \left[a_0:...:a_m
\right] \in L^m$. Let vector coordinates in $\mathbb{B}^m$ are $v = \left[v_0:...:v_n\right]$. Then
$v = L^ma$. Let $\mathfrak{M}$ be interior motion of lineal $L^m$, defined by matrix of size $(m+1)
\times(m+1)$. And let coordinates of $a$ in new basis $L'^m$ are $b = \left[b_0:...:b_m\right]$. Then
$b = \mathfrak{M}a$. Now $v = L'^mb$. Having the fact the coordinates of vector $v$ in $\mathbb{B}^n$
don't change, result matrix equality:
$$L^ma = v = L'^mb = L'^m(\mathfrak{M}a) = (L'^m\mathfrak{M})a.$$
This equality doesn't depend on vector $a$, then
\begin{equation}\label{coordchange}
L^m = L'^m\mathfrak{M}
\end{equation}
is equation of basis change in lineal.

It is necessary to find the unique form of lineal definition. Consider the following algorithm for the
unique basis search:
\begin{enumerate}
\item Let $i_p, p=\overline{0,n}$ is basis of $\mathbb{B}^n$. Start with empty basis of $L^m$.
\item Until new basis has less then $m+1$ elements, search for $i'_p$ as projection of next $i_p$ on
$L^m$.
\begin{enumerate}
\item If projection isn't null, find new vector $i''_p$ as projection of $i'_p$ on orthogonal
completion of existing basis $l_i$.
\item If $i''_p$ isn't null, find its position as free index $0 \le q \le m$ so that
$r^2 = \frac1{K_q}i''_p\odot i''_p > 0$.
\item Norm it and add to existing basis $l_q = \frac{i''_p}{r}$.
\end{enumerate}
\end{enumerate}

\section{Measure Calculus Between Lineals}
\begin{wrapfigure}{l}{0.5\textwidth}
\begin{center}
\includegraphics{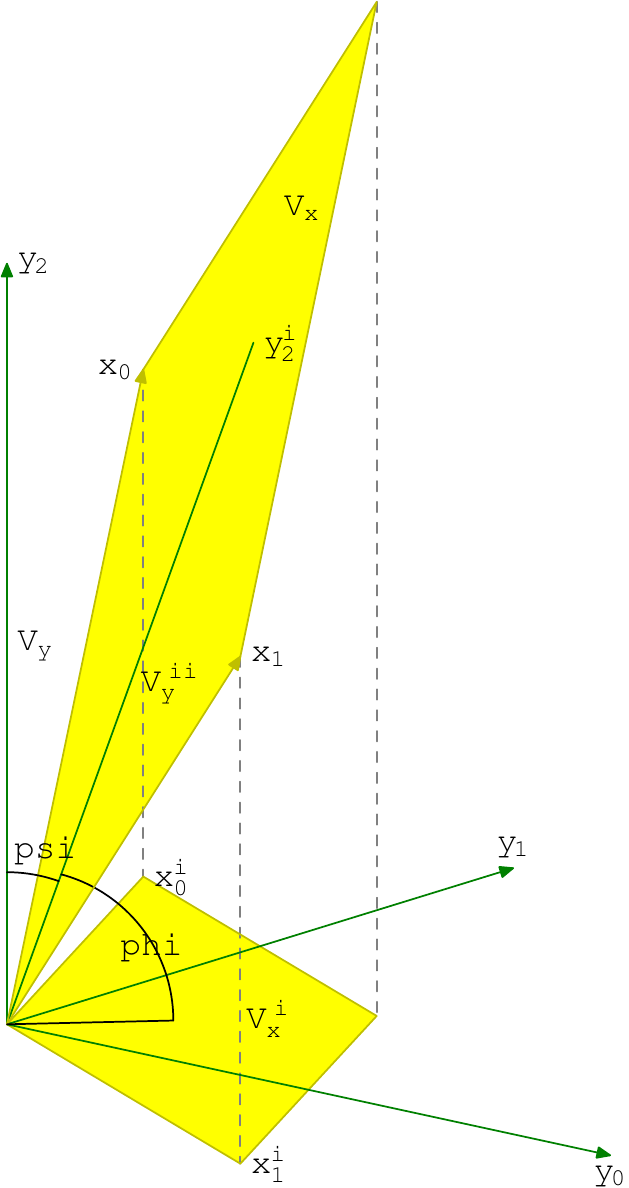}
\caption{Measure calculus between lineals $X^2$ and $Y^2$.}
\label{lineals}
\end{center}
\end{wrapfigure}
Let $X^p, Y^q, p \le q \le n$ are two lineals. Let $x_i, i = \overline{0,p}$ is the basis of $X^p$.
Let $X'^p$ be projection of $X^p$ on $Y^p$ (Figure \ref{lineals}) and let $x'_i$ be projection of $x_i$
on $Y^p$ (they are not orthonormal). If the volumes of parallelepipeds constructed on vectors $x_i$ and
$x'_i$ are equals to $V_x$ and $V'_x$ respectively and the angle between $X^p$ and $Y^q$ is measurable
and equals to $\phi$, then has place the equality:
$$V'_x = V_xC(\phi).$$
This equality is a particular case of (\ref{nearangle1}), when $k_1 = 0, T(x) = x$. In our case always
$k_1 = 0$, because the space model is linear. As were discussed earlier, $V_x = 1$ and $V'_x = \sqrt
{\det M'_x}$, where $M'_x$ is state matrix of vectors $x'_i$:
\begin{equation}\label{c_phi}
C(\phi) = \sqrt{\det M'_x}.
\end{equation}
It may happen that characteristic of $\phi$ equals to zero and we can't calculate $C^{-1}(\phi)$. If
we project vectors $x_i$ on $Y^q_\perp$ orthogonal completion of $Y^q$ (suppose it has dimension at
least $p$), get $X''^p$ constructed on vectors $x''_i$ with the volume $V''_x$, then by (\ref
{nearangle2}) get:
$$V''_x = V_xS(\phi).$$
Or, having $M''_x$ state matrix for vectors $x''_i$,
\begin{equation}\label{s_phi}
S(\phi) = \sqrt{\det M''_x}.
\end{equation}
If the dimension of $Y^q_\perp$ is less then $p$, then we can get $S(\phi)$ by projecting of $Y^q_
\perp$ on $X^p$. If $\phi$ isn't measurable, then the angle $\psi$ between $X^p$ and $Y^q_\perp$ is
measurable and:
\begin{eqnarray}
S(\psi) = \sqrt{\det M'_x},\label{s_psi}\\
C(\psi) = \sqrt{\det M''_x}.\label{c_psi}
\end{eqnarray}

The angles $\phi$ and $\psi$ present measure between lineals $X^p$ and $Y^q$. Having values $C(\phi),
S(\phi)$ and $C(\psi), S(\psi)$, it is possible to determine $\phi$ and $\psi$. The measure
characteristic of $\phi$ and $\psi$ equals to measure characteristic between $X'^p$ and $X''^p$.
Depending on this characteristic, situation can be one of the following:
\begin{itemize}
\item If characteristic equals to 1, then $\det M'_x + \det M''_x = 1$ and $\phi = \tan^{-1}\sqrt
{\frac{\det M''_x}{\det M'_x}}$, $\psi = \tan^{-1}\sqrt{\frac{\det M'_x}{\det M''_x}}$.
\item If characteristic equals to 0, then either $\det M'_x = 1$ and $\phi = \sqrt{\det M''_x}$,
$\psi = \infty$, or $\det M''_x = 1$, and $\phi = \infty$, $\psi = \sqrt{\det M'_x}$.
\item If characteristic equals to $-1$, then either $\det M'_x - \det M''_x = 1$ and $\phi = \tanh^
{-1}\sqrt{\frac{\det M''_x}{\det M'_x}}$, $\psi$ isn't measurable, or $\det M''_x - \det M'_x = 1$ and
$\phi$ isn't measurable, $\psi = \tanh^{-1}\sqrt{\frac{\det M'_x}{\det M''_x}}$, or $\det M''_x =
\det M'_x$ and $\phi = \psi = \infty$.
\end{itemize}

\section{Volume Calculation}
We can see that for any $\mathbb{B}^n$ seen as a unit sphere in $\mathbb{R}^{n+1}$, the surface is
orthogonal to radius. Let $X, Y \in \mathbb{B}^n$ and the distance between $X$ and $Y$ is small. Let
$O = (0, 0,...,0)$ is origin of $\mathbb{R}^{n+1}$. We will see that $(O-X)\odot(Y-X) = 0$ when
$Y \rightarrow X$ in sense of distance between them. $O-X = -X$, $(O-X)\odot(Y-X) = -X\odot(Y-X)
= X\odot X - X\odot Y = 1 - C_1(d(X,Y))$, where $d(X,Y)$ is distance between $X$ and $Y$. When
$Y \rightarrow X$, $d(X, Y) \rightarrow 0$ and $1 - C_1(d(X,Y)) \rightarrow 0$.

Let $A,B \in \mathbb{B}^1$. $A = \left[C_1(\alpha):S_1(\alpha)\right]$, $B = \left[C_1(\beta):
S_1(\beta)\right]$, where $C_1(x)$, $S_1(x)$ and $T_1(x)$ are defined as (\ref{cx}), (\ref{sx}) and
(\ref{tx}). Let calculate in $\mathbb{R}^2$ the area of $\mathbb{B}^1$ sector between $A$ and $B$.
In Euclidean polar system the argument $\tan\phi_e = y/x = S_1(\phi)/C_1(\phi) = T_1(\phi)$, where
$\phi$ is native argument in $\mathbb{B}^1$. The Euclidean radius $\rho = \sqrt{x^2+y^2} = \sqrt
{C_1^2(\phi) + S_1^2(\phi)}$. Having $d\phi_e = \frac{d\phi}{(1 + T_1^2(\phi))C_1^2(\phi)}$, the
area is:
$$S = \frac12\int_A^B\rho(\phi_e)^2d\phi_e
= \frac12\int_A^B(C_1^2(\phi)+S_1^2(\phi))\frac{d\phi}{C_1^2(\phi)(1+T_1^2(\phi))}$$
$$= \frac1{2}\int_A^B\frac{C_1^2(\phi)+S_1^2(\phi)}{C_1^2(\phi)+S_1^2(\phi)}d\phi
= \frac1{2}\int_A^Bd\phi = \left.\frac1{2}\phi\right|_{\alpha}^{\beta} = \frac{\beta-\alpha}{2}.$$
That is, $2S$ equals to length $AB$.

\begin{figure}[h]
\begin{center}
\includegraphics{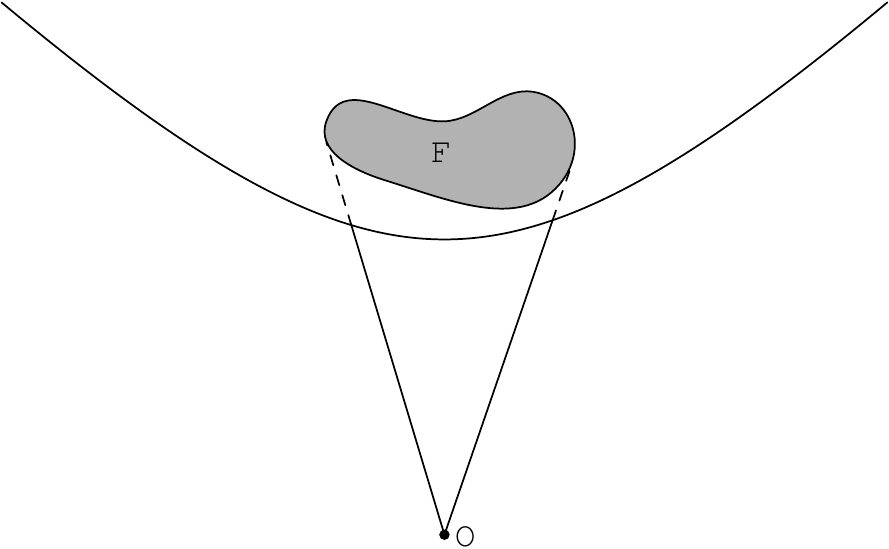}
\caption{Figure $F \subset\mathbb{B}^n$ volume calculation with aid of cone in $\mathbb{R}^{n+1}$.}
\label{cone}
\end{center}
\end{figure}

Let $F \subset \mathbb{B}^n$ be some figure with volume (in sense of $\mathbb{B}^n$) $V_{\mathbb{B}}$.
We will name $V_{\mathbb{R}}$ the volume (in sense of $\mathbb{R}^{n+1}$) of cone with base
$F \subset \mathbb{B}^n$ and vertex $O \notin \mathbb{B}^n$ origin of $\mathbb{R}^{n+1}$ (Figure
\ref{cone}). As $\mathbb{B}^n$ is orthogonal to radius, $F$ also is orthogonal. The radius equals 1,
because $\forall X \in \mathbb{B}^n, X\odot X = 1$. Then for each figure $F \subset \mathbb{B}^n$ have
place equality:

$$V_{\mathbb{B}} = (n+1)V_{\mathbb{R}}$$

As motions preserve $\mathbb{B}^n$ and the absolute value of their matrices' determinant is 1, all
motions preserve $V_{\mathbb{R}}$ and thus, they preserves also $V_{\mathbb{B}}$.

%% file: ch3-en.tex
\chapter{Theory Application}
\section{Space and Lineal Specification Search Algorithm}
As we can see, the theory described in this book is universal and easy applicable. However, one issue
stops somebody from using it. Geometric spaces are classified and defined different from the way adopted
here. Therefore, in order to not loose the feeling of reality, we will describe an algorithm aimed to
find specification for some geometric space. The algorithm can be applied to any space where have sense
notions of points, lines, planes, subspaces, distances, angles and / or motions.
\begin{enumerate}
\item Let $m$ equals to the greatest number of general situated points, or same, the lowest number of
vertices in a polyhedron of positive volume.
\item Count space dimension as $n = m - 1$.
\item Name points 0-dimensional planes and lines 1-dimensional planes.
\item For $i=\overline{1, n}$ do:
\begin{enumerate}
\item If among $(i - 1)$-dimensional planes there are non-congruent ones, then the space definition or
space terminology is inconsistent. Theory still can be used, however, in order to understand it correctly,
it is necessary to modify terminology or to define otherwise some space elements (about it later).
\item If the measure between $(i - 1)$-dimensional planes is bounded, then $k_i = 1$.
\item If the measure between $(i - 1)$-dimensional planes is scalable, then $k_i = 0$.
\item Otherwise, $k_i = -1$.
\end{enumerate}
\item Having space dimension $n$ and specification $\{k_1,..., k_n\}$, use theory.
\end{enumerate}

The necessity of proper terminology, uniform among all spaces is required by wish to have such a theory,
that isn't misleading and helps us to study the space structure and to compare it with other spaces.
Still, under inconsistent theory / terminology we should understand it has a contradiction, but failing
it to match to theory / terminology that is common today. We assume the following here:
\begin{itemize}
\item All the planes of any dimension are congruent, including points and lines.
\item Theory allows the duality principle of $(m - 1)$-dimensional planes and $(n - m)$-dimensional ones.
\end{itemize}
We should mention that `common terminology' may change over the time. In order to understand what it is
consider an example of inconsistent terminology. The Minkowskii space is successfully used in physics to
describe the theory of relativity. Unfortunately, from geometry point of view, it have no proper
terminology. The notions of `space--like lines', `time--like lines' and `light--like lines' have sense in
physics, but not in geometry. Corresponding geometric notions are: `I-st category lines', `II-nd category
lines' and `III-rd category lines'. No space motion maps some line of a category into some line of
another category. There is no contradiction here, but there is an inconsistence. What happens if somebody
wants to define a space with five categories of lines\footnote{Depending on concrete space, there can
exist more categories of two--dimensional planes. For further dimensions of planes the number of their
categories grows.}? Nobody defines several kinds of points. All points are congruent\footnote{The notion
of points on infinity is used in projective geometry. These points are non-congruent with others. The
terminology is not common in a scope of analytic geometry.}. Why shouldn't be lines all congruent? At
the other hand, \textit{relative position} of points may differ. If we name lines only the I-st category
of lines, then we should exclude II-nd and III-rd category of lines from lines. At the first look, it
conflicts with the axiom that claims any two points can be connected with a line. But this axiom may have
no place in other spaces. In contrast, just Euclidean geometry, where all points are connectable, gives
us an example of parallel lines (that have no common point). Using the duality principle, it should exist
the notion of non-connectable points (that have no common line).

It should be mention, that even for somebody feels comfortable using this theory, the algorithm described
earlier may help to determinate the specification of some exotic lineals (for example, of ones defined
as limit lineals, which aren't deductable from the space specification).

\section{Some Special Spaces}
Many linear spaces are defined using the quadric form of distance $d^2(X, Y) = (Y - X)\odot(Y - X)$. As
for these spaces $k_1 = 0$, $K_0 = 1, K_i = 0, i>0$. In this case, the equality $1 = C(d(X, Y)) = X\odot
Y = 1$ is trivial and can't be used for distance calculation. Consider one more vector product ---
$\otimes$ such as $(X\odot Y)^2 + k_1(X\otimes Y)^2 = 1, \forall X, Y\in\mathbb{B}^n, k = \{-1, 0, 1\}$.
This product is similar to exterior vector product. Change $1 = (X\odot X)(Y\odot Y)$:
\begin{eqnarray*}
(X\otimes Y)^2 &=& \frac1{k_1}((X\odot X)(Y\odot Y) - (X\odot Y)^2)\\
&=& \frac1{k_1}\left(\left(\sum_{i=0}^nK_ix_i^2\right)\left(\sum_{j=0}^nK_jy_j^2\right)
- \left(\sum_{i=0}^nK_ix_iy_i\right)\left(\sum_{j=0}^nK_jx_jy_j\right)\right)\\
&=& \frac1{k_1}\sum_{i=0}^n\sum_{j=0}^nK_iK_j(x_i^2y_j^2 - x_ix_jy_iy_j)\\
&=& \frac1{k_1}\sum_{i<j=0}^nK_iK_j(x_i^2y_j^2 + x_j^2y_i^2 - 2x_ix_jy_iy_j)\\
&=& \frac1{k_1}\sum_{i<j=0}^nK_iK_j(x_iy_j - x_jy_i)^2
\end{eqnarray*}
So,
\begin{equation}\label{vtimes}
X\otimes Y = \sqrt{\frac1{k_1}\sum_{i<j=0}^nK_iK_j(x_iy_j - x_jy_i)^2}.
\end{equation}

Note that from $\mathbb{B}^n = \{x \in \mathbb{RP}^n\,|\, x\odot x = x_0^2 = 1\}$ result $x_0 = 1$
or $x_0 = -1$. As $x \in \mathbb{B}^n$ implies $-x \in \mathbb{B}^n$ we can consider $x_0 = 1$. We
will use $\otimes$ operator for distance $d(X,Y)$ between points $X$ and $Y$. Note that having $C^2(x)
+ k_1S^2(X) = 1, \forall x\in\mathbb{R}, k_1 = \{-1, 0, 1\}$ and $(X\odot X)^2 + k_1(X\otimes Y)^2
= 1, \forall X, Y\in\mathbb{B}^n, k_1 = \{-1, 0, 1\}$ result $S(d(X, Y)) = X\otimes Y, \forall X, Y
\in\mathbb{B}^n$:
$$d^2(X,Y) = S^2(d(X,Y)) = (X\otimes Y)^2 = \frac1{k_1}\sum_{i<j=0}^nK_iK_j(x_iy_j - x_jy_i)^2$$
In this sum all non--zero terms are those for which $i = 0$:
\begin{eqnarray*}
d^2(X,Y) &=& \frac1{k_1}\sum_{j=1}^nK_j(x_0y_j - x_jy_0)^2\\
&=& \sum_{j=1}^n\frac{K_j}{k_1}(y_j - x_j)^2\\
&=& \sum_{j=1}^n(y_j - x_j)^2\prod_{p=2}^nk_p
\end{eqnarray*}
In this equality don't appear $x_0$ or $y_0$. We can consider $\mathbb{B}^n$ a hyperplane of
$\mathbb{R}^{n+1}$ with equation $x_0 = 1$ and specification $\{k_2,... k_n\}$. We can identify it with
$\mathbb{R}^n$. Then the equality above is equivalent to $(Y - X)\odot(Y - X)$.

It means firstly, that scalar product of vestors in lnear spaces ($k_1 = 0$) induces the same metrics
that is used in the model, and secondly, that non--linear spaces with specification $\{k_1, k_2,...
k_n\}$ ($k_1 \ne 0$) are best approximated by linear spaces with specification $\{0, k_2,... k_n\}$.
Note also that from here deduce that non--linear space with specification $\{k_1,... k_n\}$ is enclosed
in model meta--space of greater by one dimension, of which specification is $\{0, k_1,... k_n\}$.

We can use this quadric form to search for all characteristics except $k_1 = 0$. We will use this method
in order to describe some special spaces by specifying their specifications.

\paragraph{Case 1. Elliptic, Euclidean and Hyperbolic Spaces.}
Elliptic, linear (Euclidean) and hyperbolic (Bolyai-Lobachevsky) spaces have characteristic $k_1$ equals
to sign of space curvature $k_1 = 1$ for elliptic space, $k_1 = 0$ for linear space and $k_0 = -1$ for
hyperbolic space.

All these spaces are usually approximated by Euclidean one. We can calculate the rest of characteristics
using euclidean quadric form. Let dimension is 3:
\begin{eqnarray*}
d(X,Y)^2 &=& (y_1 - x_1)^2 + (y_2 - x_2)^2 + (y_3 - x_3)^2\\
&=& (y_1 - x_1)^2 + k_2(y_2 - x_2)^2 + k_2k_3(y_3 - x_3)^2
\end{eqnarray*}
so
$k_2 = 1$ and $k_2k_3 = 1$, $k_3 = 1$.

\paragraph{Case 2. Minkowskii Space.}
The distance between $X$ and $Y$ is calculated (for time--like vectors) as
\begin{eqnarray*}
d^2(X,Y) &=& (y_1 - x_1)^2 - (y_2 - x_2)^2 - (y_3 - x_3)^2 - (y_4 - x_4)^2\\
&=& (y_1 - x_1)^2 + k_2(y_2 - x_2)^2 + k_2k_3(y_3 - x_3)^2 + k_2k_3k_4(y_4 - x_4)^2
\end{eqnarray*}
where coordinate 1 is time--like and coordinates 2, 3 and 4 are space--like. So $k_2 = -1$, $k_2k_3 =
-1$, $k_3 = 1$ and $k_2k_3k_4 = -1$, $k_4 = 1$. As Minkowskii space is linear, $k_1 = 0$.

If we introduce curvature in space, its structure changes. For example, let $k_1 = 1$. Then
\begin{eqnarray*}
X\odot Y &=& x_0y_0 + k_1x_1y_1 + k_1k_2x_2y_2 + k_1k_2k_3x_3y_3 + k_1k_2k_3k_4x_4y_4\\
&=& x_0y_0 + x_1y_1 - x_2y_2 - x_3y_3 - x_4y_4
\end{eqnarray*}
so time characteristic becomes elliptic and space characteristic becomes hyperbolic. If $k_1 = 1$, then
\begin{eqnarray*}
X\odot Y &=& x_0y_0 + k_1x_1y_1 + k_1k_2x_2y_2 + k_1k_2k_3x_3y_3 + k_1k_2k_3k_4x_4y_4\\
&=& x_0y_0 - x_1y_1 + x_2y_2 + x_3y_3 + x_4y_4
\end{eqnarray*}
and time characteristic becomes hyperbolic and space characteristic becomes elliptic.

\paragraph{Case 3. Minkowskii Space with 2-dimensional Time.}
Consider a 4-dimensional space with distance quadric form that has 2 positive signs and 2 negative.
This space is sometimes named Minkowskii space with 2-dimensional time:
\begin{eqnarray*}
d^2(X,Y) &=& (y_1 - x_1)^2 + (y_2 - x_2)^2 - (y_3 - x_3)^2 - (y_4 - x_4)^2\\
&=& (y_1 - x_1)^2 + k_2(y_2 - x_2)^2 + k_2k_3(y_3 - x_3)^2 + k_2k_3k_4(y_4 - x_4)^2
\end{eqnarray*}
So $k_2 = 1$, $k_2k_3 = -1$, $k_3 = -1$ and $k_2k_3k_4 = -1$, $k_4 = 1$. As for all lineal spaces, $k_1
= 0$ for it.

\paragraph{Case 4. Spaces with Degenerate Distance Quadric Form.}
Consider linear 4-dimensional space ($k_1 = 0$) with degenerate distance quadric form:
\begin{eqnarray*}
d^2(X,Y) &=& (y_1 - x_1)^2 + (y_2 - x_2)^2 + (y_3 - x_3)^2\\
&=& (y_1 - x_1)^2 + k_2(y_2 - x_2)^2 + k_2k_3(y_3 - x_3)^2 + k_2k_3k_4(y_4 - x_4)^2
\end{eqnarray*}
so $k_2 = k_3 = 1$ and $k_4 = 0$.

It means that motions:
\begin{eqnarray*}
x'_1 &=& x_1\\
x'_2 &=& x_2\\
x'_3 &=& x_3\\
x'_4 &=& \phi_1x_1 + \phi_2x_2 + \phi_3x_3 + x_4
\end{eqnarray*}
are all valid.

However, transformation:
\begin{eqnarray*}
x'_1 &=& x_1\\
x'_2 &=& x_2\\
x'_3 &=& x_3\\
x'_4 &=& \phi_4x_4
\end{eqnarray*}
is not a motion. Although it preserved distance, it doesn't preserve volume except $\phi_4 = 1$ or
$-1$. It is an example of angle scaling.

\section{Spaces as Product of their Subspaces}
Another way to define spaces is by product of their subspaces. It is necessary to be accurate here. The
geometric space isn't only a structure of points. It is also the structure of all its subspaces. It is
mistake to think that having $\mathbb{R}^1$ is isomorphic to one--dimensional Euclidean space
$\mathbb{E}^1$, from $\mathbb{R}^1\times\mathbb{R}^1 = \mathbb{R}^2$ results $\mathbb{E}^1\times
\mathbb{E}^1 = \mathbb{E}^2$ (using specification notation, $\{0\}\times$\{0\} = \{0, 1\}). The
problem is the product doesn't define way to measure the angle between multiplied subspaces. It can be
defined in several ways, for example, $\{0, 0\}$ or $\{0, -1\}$.

The situation is even worse when multiplied subspaces with different specification $X^m$ and $Y^n$.
One--dimensional images can be constructed in two ways: $X^1\times Y^0$ (isomorphic to $X^1$) and $X^0
\times Y^1$ (isomorphic to $Y^1$). And if $X^1$ and $Y^1$ have different specifications, these two
one--dimensional lines aren't congruent. For example, if somebody wants to construct geometry on a
cylinder, first thing he or she thinks of is $\mathbb{S}^1\times\mathbb{E}^1$ ($\{1\}\times\{0\}$),
where $\mathbb{S}^1$ is one--dimensional elliptic space. In this case some lines are circles, some are
lines and others are right and left helices, that may not intersect, intersect in one point or intersect
in infinity of points. As an example of complete geometry on a cylinder you can take the space with
specification $\{1, 0\}$.

Additionally, one should not consider that if from algebraic point of view $\mathbb{E}^1$ is isomorphic
to $\mathbb{H}^1$ (one--dimensional hyperbolic space), then constructions like $\mathbb{H}^2\times
\mathbb{E}^1$ ($\{-1, 1\}\times\{0\}$) and $\mathbb{H}^2\times \mathbb{H}^1$ ($\{-1, 1\}\times
\{-1\}$) are also isomorphic. From geometric point of view, $\mathbb{E}^1$ is scalable, while $\mathbb
{H}^1$ is not (the mutual departure of points is possible, however it can't be linear). In contrast, it
is possible to construct spaces with specifications $\{-1, 1, 0\}$ and $\{-1, 1, -1\}$, which differ one
from another by the fact that in first one on a twodimensional plane doesn't containing some point there
is the only point that isn't connectable with it, and for the second space the number of such a points is
infinity.